\documentclass[10pt]{amsart}
\usepackage[margin=2.5cm]{geometry}
\usepackage{tikz-cd}
\usepackage{graphicx}					
\usepackage{amssymb}
\usepackage{amsmath}
\usepackage{subfig, graphicx}
\usepackage{mathrsfs}
\usepackage{amssymb, amsthm, indentfirst}
\usepackage{relsize}
\usepackage{hyperref}
\usepackage{stmaryrd}
\usepackage{lineno}

\numberwithin{equation}{section}
\usepackage{color}
\theoremstyle{plain}
\newtheorem{thm}{Theorem}[section]

\newtheorem{conj}[thm]{Conjecture}
\newtheorem{lemma}[thm]{Lemma}

\theoremstyle{definition}

\newtheorem{rmk}[thm]{Remark}

\newtheorem*{hypothesis*}{Hypothesis}

\author{SHIH-YU CHEN}
\address{Institute of Mathematics~\\Academia Sinica~\\ 6F, Astronomy-Mathematics Building, No.\,1, Sec.\,4, Roosevelt Road, Taipei 10617, Taiwan, ROC}
\email{sychen0626@gate.sinica.edu.tw}

\def\GL{{\rm{GL}}}

\def\o{\frak{o}}

\def\A{{\mathbb A}}
\def\C{{\mathbb C}}
\def\F{{\mathbb F}}
\def\R{{\mathbb R}}
\def\Q{{\mathbb Q}}
\def\Z{{\mathbb Z}}

\def\<{\langle}
\def\>{\rangle}

\def\bp{\begin{pmatrix}}
\def\ep{\end{pmatrix}}

\def\<{\langle}
\def\>{\rangle}

\def\GL{\operatorname{GL}}

\def\1{\mathbf{1}}

\def\itPi{\mathit{\Pi}}

\def\itSigma{\mathit{\Sigma}}

\title{Algebraicity of the near central non-critical values of symmetric fourth $L$-functions for Hilbert modular forms}

\begin{document}
\begin{abstract}
Let $\itPi$ be a cohomological irreducible cuspidal automorphic representation of $\GL_2(\A_\F)$ with central character $\omega_\itPi$ over a totally real number field $\F$. In this paper, we prove the algebraicity of the near central non-critical value of the symmetric fourth $L$-function of $\itPi$ twisted by $\omega_\itPi^{-2}$. The algebraicity is expressed in terms of the Petersson norm of the normalized newform of $\itPi$ and the top degree Whittaker period of the Gelbart--Jacquet lift ${\rm Sym}^2\itPi$ of $\itPi$.
\end{abstract}

\maketitle
\tableofcontents

\section{Introduction}
Let $\F$ be a totally real number field with $[\F:\Q]=d$.
Let $\itPi = \bigotimes_v \itPi_v$ be a cohomological irreducible cuspidal automorphic representation of $\GL_2(\A_\F)$ with central character $\omega_\itPi$. We have $|\omega_\itPi| = |\mbox{ }|_{\A_\F}^{{\sf w}}$ for some ${\sf w} \in \Z$.
Denote by $\Q(\itPi)$ the Hecke field of $\itPi$. 
For each archimedean place $v$ of $\F$, we have
\[
\itPi_v = D_{\kappa_v} \otimes |\mbox{ }|^{{\sf w}/2}
\]
for some $\kappa_v \in \Z_{\geq 2}$ such that $\kappa_v \equiv {\sf w} \,({\rm mod}\,2)$, where $D_{\kappa_v}$ is the discrete series representation of $\GL_2(\R)$ of weight $\kappa_v$.
Let $f_\itPi$ and $f_{\itPi^\vee}$ be the normalized newforms of $\itPi$ and $\itPi^\vee$, respectively.
Let $\Vert f_\itPi \Vert$ be the Petersson norm of $f_\itPi$ defined by
\begin{align*}
\Vert f_\itPi \Vert = \int_{\A_\F^\times\GL_2(\F)\backslash \GL_2(\A_\F)}f_\itPi(g)f_{\itPi^\vee}(g\cdot {\rm diag}(-1,1)_\infty)\,dg^{\rm Tam},
\end{align*}
where $dg^{\rm Tam}$ is the Tamagawa measure.
Let ${\rm Sym}^2\itPi$ be the functorial lift of $\itPi$ with respect to the symmetric square representation of $\GL_2$. The functoriality was established by Gelbart--Jacquet \cite{GJ1978}.
Suppose that $\itPi$ is non-CM. Then ${\rm Sym}^2\itPi$ is a  cohomological irreducible cuspidal automorphic representation of $\GL_3(\A_\F)$.
We have the $\Q(\itPi)$-rational structure on ${\rm Sym}^2\itPi_f = \bigotimes_{v \nmid \infty}{\rm Sym}^2\itPi_v$ defined through the Whittaker model of ${\rm Sym}^2\itPi$ by fixing specific normalization of archimedean Whittaker functions in the minimal ${\rm SO}(3)$-type of ${\rm Sym}^2\itPi_v$ following Miyazaki \cite{Miyazaki2009} (cf.\,Theorem \ref{T:Miyazaki}) for each archimedean place $v$ of $\F$. 
On the other hand, ${\rm Sym}^2\itPi$ contributes to the cuspidal cohomology of $\GL_3(\A_\F)$ with coefficients in certain irreducible algebraic representation of ${\rm R}_{\F/\Q}\GL_{3/\F}$. 
This in term defines $\Q(\itPi)$-rational structures on ${\rm Sym}^2\itPi_f$ in the bottom and top degree of cuspidal cohomology. 
Comparing the $\Q(\itPi)$-rational structures, we obtain non-zero complex numbers $p^b({\rm Sym}^2\itPi)$ and $p^t({\rm Sym}^2\itPi)$, well-defined up to $\Q(\itPi)^\times$, called the bottom degree and top degree Whittaker periods of ${\rm Sym}^2\itPi$, respectively.
It is known that the bottom degree Whittaker period is related to the algebraicity of the critical values for $\GL_3$ and $\GL_3 \times \GL_2$ by Mahnkopf \cite{Mahnkopf2005} and Raghuram \cite{Raghuram2009} (see also \cite{RS2017}, \cite{Sachdeva2020}, and Conjecture \ref{C:Deligne-Whittaker} below).
In this paper, we prove similar result for the top degree Whittaker period and its connection to non-critical values.
More precisely, let 
\[
L(s, \itPi, {\rm Sym}^4\otimes \omega_\itPi^{-2})
\]
be the twisted symmetric fourth $L$-function of ${\rm Sym}^4\itPi$ by $\omega_\itPi^{-2}$. We denote by $L^{(\infty)}(s, \itPi,{\rm Sym}^4 \otimes \omega_\itPi^{-2})$ the $L$-function obtained by excluding the archimedean $L$-factors.
The following theorem is our main result which expresses the algebraicity of the non-critical value $L^{(\infty)}(1,\itPi,{\rm Sym}^4 \otimes \omega_\itPi^{-2})$ in terms of the top degree Whittaker period of ${\rm Sym}^2\itPi$ and the Petersson norm $\Vert f_\itPi \Vert$.

\begin{thm}\label{T:main}
Suppose $\itPi$ is non-CM and $\kappa_v \geq 3$ for all archimedean places $v$ of $\F$. We have
\begin{align*}
\sigma \left( \frac{L^{(\infty)}(1,\itPi,{\rm Sym}^4 \otimes \omega_\itPi^{-2})}{\pi^{3\sum_{v \mid \infty}\kappa_v}\cdot G(\omega_\itPi)^{-3}\cdot \Vert f_\itPi\Vert\cdot p^t({\rm Sym}^2\itPi)}\right) = \frac{L^{(\infty)}(1,{}^\sigma\!\itPi,{\rm Sym}^4 \otimes {}^\sigma\!\omega_\itPi^{-2})}{\pi^{3\sum_{v \mid \infty}\kappa_v}\cdot G({}^\sigma\!\omega_\itPi)^{-3}\cdot \Vert f_{{}^\sigma\!\itPi}\Vert\cdot p^t({\rm Sym}^2{}^\sigma\!\itPi)}
\end{align*}
for all $\sigma \in {\rm Aut}(\C)$. Here $G(\omega_\itPi)$ is the Gauss sum of $\omega_\itPi$ and $p^t({\rm Sym}^2\itPi)$ is the top degree Whittaker period of ${\rm Sym}^2\itPi$.
In particular, we have
\[
\frac{L^{(\infty)}(1,\itPi,{\rm Sym}^4 \otimes \omega_\itPi^{-2})}{\pi^{3\sum_{v \mid \infty}\kappa_v}\cdot G(\omega_\itPi)^{-3}\cdot \Vert f_\itPi\Vert\cdot p^t({\rm Sym}^2\itPi)} \in \Q(\itPi).
\]
\end{thm}

\begin{rmk}
We have similar results \cite{Urban1995}, \cite{Hida1999}, and \cite{Namikawa2015} for symmetric square $L$-functions for $\GL_2$ over number fields, in which cases the algebraicity are expressed in terms of the product of the bottom and top degree Whittaker periods for $\GL_2$. 
These results are vastly generalized by Balasubramanyam--Raghuram in \cite{BR2017} to $\GL_n$. 
A crude form of Theorem \ref{T:main} is also given in \cite[(3.4.3)]{BR2017}.
See also the results of Grobner--Harris--Lapid \cite[\S\,6]{GHL2016} and Grobner \cite{Grobner2018} on the algebraicity of non-critical $L$-values in terms of top degree Whittaker and Shalika periods.
\end{rmk}

\subsection{A sketch of the proof}

There are two main ingredients for the proof of Theorem \ref{T:main}:
\begin{itemize}
\item algebraicity of the Rankin--Selberg $L$-functions for $\GL_3 \times \GL_2$;
\item algebraicity of the adjoint $L$-functions for $\GL_3$.
\end{itemize}
Let $\itSigma$ be a cohomological irreducible cuspidal automorphic representation of $\GL_3(\A_\F)$. 
By comparing the $\Q(\itSigma)$-rational structures on $\itSigma_f = \bigotimes_{v \nmid \infty}\itSigma_v$ given by the Whittaker model and by the cuspidal cohomology for $\GL_3$ in either bottom degree of top degree, we obtain the bottom degree and top degree Whittaker periods 
\[
p^b(\itSigma) \in \C^\times,\quad p^t(\itSigma) \in \C^\times,
\]
which are uniquely determined up to $\Q(\itSigma)^\times$.
Similarly we have the Whittaker periods
\[
p(\itPi,\underline{\varepsilon}) \in \C^\times
\]
of $\itPi$ which are indexed by sequence of signs $\underline{\varepsilon} = (\varepsilon_v)_{v \nmid \infty}$.
The definition of the Whittaker periods depends on the choice of archimedean Whittaker functions of $\itPi_\infty = \bigotimes_{v \mid \infty}\itPi_v$ and $\itSigma_\infty = \bigotimes_{v \mid \infty}\itSigma_v$ in the minimal ${\rm O}(2)$-type and the minimal ${\rm O}(3)$-type, respectively.
For $\itPi_\infty$, we follow the classical normalization (cf.\,(\ref{E:GL_2 Whittaker})). 
As for $\itSigma_\infty$,we follow the normalization due to Miyazaki \cite{Miyazaki2009} (cf.\,Theorem \ref{T:Miyazaki} and (\ref{E:GL_3 class})).
Under certain balanced condition on the cuspidal parameters of $\itPi$ and $\itSigma$, Raghuram \cite{Raghuram2016} proved the algebraicity of the critical values of the Rankin--Selberg $L$-function $L(s,\itSigma \times \itPi)$ in terms of the Whittaker periods $p^b(\itSigma)$ and $p(\itPi,\underline{\varepsilon})$. In \cite{BR2017}, Balasubramanyam--Raghuram proved the algebraicity of the residue ${\rm Res}_{s=1}L(s,\itSigma\times\itSigma^\vee)$ of the Rankin--Selberg $L$-function of $\itSigma \times \itSigma^\vee$ in terms of the Whittaker periods $p^t(\itSigma)$ and $p^b(\itSigma^\vee)$. More precisely, there exist non-zero complex numbers $p(m,\itSigma_\infty \times \itPi_\infty)$ for each critical point $m+\tfrac{1}{2} \in \Z+\tfrac{1}{2}$ of $L(s,\itSigma \times \itPi)$ and $p(\itSigma_\infty,{\rm Ad})$ such that
\begin{align}
\frac{L^{(\infty)}(m+\tfrac{1}{2},\itSigma \times \itPi)}{|D_\F|^{1/2}\cdot p(m,\itSigma_\infty \times \itPi_\infty)\cdot G(\omega_\itPi)\cdot p^b(\itSigma)\cdot p(\itPi,\underline{\varepsilon}(m,\itSigma))} \in \Q(\itSigma)\cdot\Q(\itPi),\label{E:1.1}\\
\frac{{\rm Res}_{s=1}L^{(\infty)}(s,\itSigma\times\itSigma^\vee)}{{\rm Reg}_\F\cdot p(\itSigma_\infty,{\rm Ad})\cdot p^t(\itSigma)\cdot p^b(\itSigma^\vee)} \in \Q(\itSigma)\cdot\Q(\itPi).\label{E:1.2}
\end{align}
Here $\underline{\varepsilon}(m,\itSigma)$ is certain sequence of signs uniquely determined by $m$ and $\itSigma$.
The archimedean factors $p(m,\itSigma_\infty \times \itPi_\infty)$ and $p(\itSigma_\infty,{\rm Ad})$ are defined to be the inverse of certain weighted sum of archimedean local zeta integrals involving Whittaker functions of $\itSigma_\infty$ and $\itPi_\infty$. In Theorems \ref{T:RS} and \ref{T:adjoint} below, we determine the rationality of these archimedean factors up to powers of $2\pi\sqrt{-1}$, which is the main novelty of this paper.
Now we take $\itSigma = {\rm Sym}^2\itPi$ be the Gelbart--Jacquet lift of $\itPi$, which is cuspidal by our assumption that $\itPi$ is non-CM. 
Consider the triple product $L$-function $L(s,\itPi \times \itPi \times \itPi)$, which has the factorization
\[
L(s,\itPi \times \itPi \times \itPi) = L(s, {\rm Sym}^2\itPi \times \itPi)\cdot L(s, \itPi \times \omega_\itPi).
\]
We have the results of Garrett--Harris \cite{GH1993} and Raghuram--Tanabe \cite{RT2011} on the algebraicity of the triple product $L$-functions in the balanced case and the twisted standard $L$-functions for $\GL_2$, respectively. Comparing with (\ref{E:1.1}), together with the assumption $\kappa_v \geq 3$ for all archimedean places $v$, we obtain the period relation (cf.\,Theorem \ref{T:GJ period relation})
\begin{align}\label{E:1.3}
\frac{|D_\F|^{1/2}\cdot\pi^{2d}\cdot G(\omega_\itPi)^{-3}\cdot\Vert f_\itPi \Vert^2}{p^b({\rm Sym}^2\itPi^\vee)} \in \Q(\itPi).
\end{align}
Finally, note that
\[
L(s,{\rm Sym}^2\itPi \times {\rm Sym}^2\itPi^\vee) = L(s,\itPi,{\rm Sym}^4 \otimes \omega_\itPi^{-2})\cdot L(s,\itPi ,{\rm Sym}^2\otimes \omega_\itPi^{-1})\cdot \zeta_\F(s).
\]
Theorem \ref{T:main} then follows from (\ref{E:1.2}), (\ref{E:1.3}), Deligne's conjecture for $L(1,\itPi ,{\rm Sym}^2\otimes \omega_\itPi^{-1})$, and the class number formula.

This paper is organized as follows. In \S\,\ref{S:Whittaker period}, we recall the definition of Whittaker periods. The periods are obtained by comparing rational structures via the Whittaker model and via the cuspidal cohomology. The key point here is that in the definition of the classes (\ref{E:GL_3 class}) in the bottom and top degree cuspidal cohomology groups, we follow the normalization of archimedean Whittaker functions by Miyazaki \cite{Miyazaki2009}. In \S\,\ref{S:proof}, we prove our main result based on Theorems \ref{T:GJ period relation} and \ref{T:adjoint}. In \S\,\ref{S:GL_3 GL_2}, we study the algebraicity of Rankin--Selberg $L$-functions for $\GL_3 \times \GL_2$. The main result of \S\,\ref{S:GL_3 GL_2} is Theorem \ref{T:RS}, which refines the result of Raghuram \cite{Raghuram2009}. As a consequence, we obtain two period relations in Theorem \ref{T:GJ period relation} and Conjecture \ref{C:Deligne-Whittaker}. One is a period relation between the bottom degree Whittaker period of the Gelbart--Jacquet liftings and the Petersson norm of the normalized newform for $\GL_2$. The other is a conjectural period relation between the bottom degree Whittaker periods for $\GL_3$ and the product of Deligne's motivic periods for the conjectural motive attached to cohomological irreducible cuspidal automorphic representations of $\GL_3$.
In \S\,\ref{S:adjoint}, we study the algebraicity of the adjoint $L$-values at $s=1$ for $\GL_3$. We prove in Theorem \ref{T:adjoint} a refinement of the result of Balasubramanyam--Raghuram \cite{BR2017}.

\subsection{Notation}
Let $\GL_n$ be the general linear group of rank $n$ over $\Q$. Let $N_n$ and $T_n$ be the maximal unipotent subgroup and the maximal torus of $\GL_n$ consisting of upper triangular unipotent matrices and diagonal matrices, respectively.
For $a_1,\cdots,a_n \in \GL_1$, we write
\[
{\rm diag}(a_1,\cdots,a_n) = \bp a_1 &\cdots &0\\ \vdots&\ddots&\vdots \\ 0&\cdots& a_n \ep \in T_n.
\]
Denote by $X^+(T_n)$ the set of dominant integral weights of $T_n$ with respect to the standard Borel subgroup consisting of upper triangular matrices. We regard $X^+(T_n)$ as the set of $n$-tuples of integers $\mu = (\mu_1,\cdots,\mu_n)$ such that $\mu_1 \geq \cdots \geq \mu_n$. For $\mu \in X^+(T_n)$, its dual weight $\mu^\vee \in X^+(T_n)$ is defined by $\mu^\vee_i = -\mu_{n-i+1}$ for $1 \leq i \leq n$.

Let $K_n$ be the compact modulo center subgroup of $\GL_n(\R)$ defined by
\[
K_{n} = \R_+\cdot {\rm SO}(n).
\]
Here we regard the set $\R_+$ of positive real numbers as the topological connected component of the center of $\GL_n(\R)$. 
We denote by $\frak{g}_n$ and $\frak{k}_n$ the Lie algebras of $\GL_n(\R)$ and $K_n$, respectively. Their complexifications are denoted by $\frak{g}_{n,\C}$ and $\frak{k}_{n,\C}$. Let $\frak{g}_{n,\C}^*$ and $\frak{k}_{n,\C}^*$ be the dual spaces of $\frak{g}_{n,\C}$ and $\frak{k}_{n,\C}$, respectively. 
For $n=2$, let $Y_\pm \in \frak{g}_{2,\C} / \frak{k}_{2,\C}$ defined by
\[
Y_\pm = \bp 1 & 0\\ 0& -1\ep \otimes \sqrt{-1} \pm \bp 0 & -1 \\ -1 & 0 \ep,
\]
and $\{Y_+^*,Y_-^*\} \subset (\frak{g}_{2,\C} / \frak{k}_{2,\C})^*$ be the corresponding dual basis.
For $n=3$, let $\{X_2,X_1,X_0,X_{-1},X_{-2}\}$ be the basis of $\frak{g}_{n,\C} / \frak{k}_{n,\C}$ defined by
\begin{align}\label{E:Lie algebra basis}
\begin{split}
X_{\pm2} &= \bp 1 & 0 & 0\\ 0& -1 & 0 \\ 0&0&0\ep \otimes \sqrt{-1} \pm \bp 0 & -1 & 0 \\ -1 & 0 & 0 \\ 0&0&0 \ep,\quad 
X_0 = \frac{1}{3} \bp -1&0&0 \\ 0&-1&0 \\ 0&0&2 \ep\otimes \sqrt{-1},\\ 
X_{\pm1} & = \pm\frac{1}{2} \bp 0&0&1 \\ 0&0&0 \\ 1&0&0 \ep \otimes \sqrt{-1} + \frac{1}{2}\bp 0&0&0\\0&0&-1\\0&-1&0\ep,
\end{split}
\end{align}
and $\{X_2^*,X_1^*,X_0^*,X_{-1}^*,X_{-2}^*\} \subset (\frak{g}_{n,\C} / \frak{k}_{n,\C})^*$ be the corresponding dual basis.

Let $\F$ be a totally real number field. 
Let $S_\infty$ be the set of archimedean places of $\F$.
Let $D_\F$ and ${\rm Reg}_\F$ be the discriminant and regulator of $\F$, respectively.
For each place $v$ of $\F$, let $\F_v$ be the completion of $\F$ at $v$ and $\frak{o}_{\F_v}$ the ring of integers of $\F_v$ when $v$ is finite.
For $n \in \Z_{\geq 1}$, let
\[
\frak{g}_{n,\infty} = \bigoplus_{v \in S_\infty} \frak{g}_n,\quad \frak{k}_{n,\infty} = \bigoplus_{v \in S_\infty} \frak{k}_n
\]
be the Lie algebras of the real Lie groups
\[
\GL_n(\F_\infty) = \prod_{v \in S_\infty}\GL_n(\F_v),\quad K_{n,\infty} = \prod_{v \in S_\infty} K_n,
\]
respectively.
Let $\A_\F$ be the ring of adeles of $\F$ and $\A_{\F,f}$ its finite part. 
Let $\psi_\Q = \bigotimes_v \psi_{\Q_v} : \Q\backslash\A_\Q \rightarrow \C$ be the additive character defined so that
\begin{align*}
\psi_{\Q_p}(x) & = e^{-2\pi \sqrt{-1}x} \mbox{ for }x \in \Z[p^{-1}],\\
\psi_{\R}(x) & = e^{2\pi \sqrt{-1}x} \mbox{ for }x \in \R.
\end{align*}
Let $\psi_{n,\F}=\bigotimes_v \psi_{n,\F_v} : N_n(\F)\backslash N_n(\A_\F) \rightarrow \C$ be the non-degenerated additive character defined by
\[
\psi_{n,\F}(u) = \psi_\Q\circ {\rm tr}_{\F/\Q}(u_{12}+u_{23}+\cdots u_{n-1,n})
\]
for $u=(u_{ij}) \in N_n(\A_\F)$.
We write $\psi_\F = \psi_{1,\F}$.
Let 
\[
\Gamma_\R(s) = \pi^{-s/2}\Gamma\left(\tfrac{s}{2}\right),\quad \Gamma_\C(s) = 2(2\pi)^{-s}\Gamma(s),
\]
where $\Gamma(s)$ is the gamma function.
We recall Barne's first and second lemmas, which will be used in the proof of Lemmas \ref{L:archimedean RS} and \ref{L:archimedean adjoint} below.
\begin{lemma}
Let $a,b,c,d,e \in \C$.
\begin{itemize}
\item[(1)] We have
\begin{align*}
\int_L \frac{ds}{2\pi\sqrt{-1}}\,\Gamma_\R(s+a)\Gamma_\R(s+b)\Gamma_\R(-s+c)\Gamma_\R(-s+d) & = 2 \cdot \frac{\Gamma_\R(a+c)\Gamma_\R(a+d)\Gamma_\R(b+c)\Gamma_\R(b+d)}{\Gamma_\R(a+b+c+d)}.
\end{align*}
Here $L$ is any vertical path from south to north which keep the poles of $\Gamma_\R(s+a)\Gamma_\R(s+b)$ and $\Gamma_\R(-s+c)\Gamma_\R(-s+d)$ on its left and right, respectively.
\item[(2)] We have
\begin{align*}
&\int_L \frac{ds}{2\pi\sqrt{-1}}\,\frac{\Gamma_\R(s+a)\Gamma_\R(s+b)\Gamma_\R(s+c)\Gamma_\R(-s+d)\Gamma_\R(-s+e)}{\Gamma_\R(s+a+b+c+d+e)} \\
& =2\cdot \frac{\Gamma_\R(a+d)\Gamma_\R(a+e)\Gamma_\R(b+d)\Gamma_\R(b+e)\Gamma_\R(c+d)\Gamma_\R(c+e)}
{\Gamma_\R(b+c+d+e)\Gamma_\R(a+c+d+e)\Gamma_\R(a+b+d+e)}
\end{align*}
Here $L$ is any vertical path from south to north which keep the poles of $\Gamma_\R(s+a)\Gamma_\R(s+b)\Gamma_\R(s+c)$ and $\Gamma_\R(-s+d)\Gamma_\R(-s+e)$ on its left and right, respectively.
\end{itemize}
\end{lemma}

\subsection{Measures}\label{SS:measure}
Let $v$ be a place of $\F$. The Haar measure $dx_v$ on $\F_v$ is normalized so that ${\rm vol}(\o_{\F_v},dx_v)=1$ if $v$ is finite, and $dx_v$ is the Lebesgue measure on $\F_v$ if $\F_v = \R$. The Haar measure $d^\times x_v$ on $\F_v^\times$ is normalized so that ${\rm vol}(\o_{\F_v}^\times,d^\times x_v)=1$ if $v$ is finite, and 
\[
d^\times x_v =\frac{dx_v}{|x_v|}
\]
if $v$ is archimedean.
The Haar measure $dg = \prod_{v}dg_v$ on $\GL_n(\A_\F)$ is defined as follows: When $v$ is finite, we assume ${\rm vol}(\GL_n(\o_{\F_v}),dg_v)=1$. When $v$ is archimedean, we have
\begin{align*}
\int_{\GL_n(\F_v)}f(g_v)\,dg_v = \int_{N_n(\F_v)}du\int_{\R_{+}^n}d^\times a_1\cdots d^\times a_n\int_{{\rm O}(n)}dk\,f(u\cdot{\rm diag}(a_1\cdots a_n,a_2\cdots a_n,\cdots,a_n)\cdot k)
\end{align*}
for $f \in L^1(\GL_n(\F_v))$. Here $du$ is defined by the product measure on $\F_v^{n(n-1)/2}$ and ${\rm vol}({\rm O}(n),dk)=1$.

\subsection{Gauss sums}
Let $\chi$ be an algebraic Hecke character of $\A_\F^\times$.
The signature ${\rm sgn}(\chi)$ of $\chi$ is the sequence 
\[
{\rm sgn}(\chi) = (\chi_v(-1))_{v \in S_\infty} \in \{\pm1\}^{S_\infty}.
\]
We say that $\chi$ has parallel signature if $\chi_{v_1}(-1) = \chi_{v_2}(-1)$ for all $v_1,v_2 \in S_\infty$.
The Gauss sum $G(\chi)$ of $\chi$ is defined by
\[
G(\chi) = |D_\F|^{-1/2}\prod_{v \nmid \infty}\varepsilon(0,\chi_v,\psi_{\F_v}),
\]
where $\varepsilon(s,\chi_v,\psi_{\F_v})$ is the $\varepsilon$-factor of $\chi_v$ with respect to $\psi_{\F_v}$ defined in \cite{Tate1979}.
For $\sigma \in {\rm Aut}(\C)$, define Hecke character ${}^\sigma\!\chi$ of $\A_\F^\times$ by ${}^\sigma\!\chi(x) = \sigma(\chi(x))$.
It is easy to verify that 
\begin{align}\label{E:Galois Gauss sum}
\begin{split}
\sigma(G(\chi)) &= {}^\sigma\!\chi(u_\sigma)G({}^\sigma\!\chi),\\
\sigma\left(\frac{G(\chi\chi')}{G(\chi)G(\chi')}\right) &= \frac{G({}^\sigma\!\chi{}^\sigma\!\chi')}{G({}^\sigma\!\chi)G({}^\sigma\!\chi')}
\end{split}
\end{align}
for algebraic Hecke characters $\chi,\chi'$ of $\A_\F^\times$,
where $u_\sigma \in \prod_p \Z_p^\times \subset \A_{\F,f}^\times$ is the unique element such that $\sigma(\psi_{\Q}(x)) = \psi_{\Q}(u_\sigma x)$ for $x \in \A_{\Q,f}$.

\section{Whittaker periods for $\GL_3$}\label{S:Whittaker period}

\subsection{Irreducible algebraic representations}

\subsubsection{Irreducible algebraic representations of ${\rm R}_{\F/\Q}\GL_{n/\F}$}

For $\mu\in X^+(T_n)$, let $(\rho_\mu,M_\mu)$ be the irreducible algebraic representation of $\GL_n$ with highest weight $\mu$. Let $M_{\mu,\C} = M_\mu \otimes_\Q \C$ be the base change to a representation of $\GL_n(\C)$. 
Note that $M_\mu^\vee = M_{\mu^\vee}$.

Let
\[
\mu = \prod_{v \in S_\infty} \mu_v \in \prod_{v \in S_\infty}X^+(T_n)
\]
be a sequence of dominant integral weights of $T_n$ indexed by $S_\infty$. Write $\mu_v = (\mu_{1,v},\cdots,\mu_{n,v})$.  We say $\mu$ is pure if there exists an integer ${\sf w}$ such that $\mu_{i,v}+\mu_{n-i+1,v} = {\sf w}$ for all $1 \leq i \leq n$ and $v \in S_\infty$.
Let 
\[
(\rho_\mu,M_{\mu,\C}) = (\otimes_{v \in S_\infty}\rho_{\mu_v},\otimes_{v \in S_\infty}M_{\mu_v,\C}).
\]  
For $\sigma \in {\rm Aut}(\C)$, let ${}^\sigma\!\mu = \prod_{v \in S_\infty}{}^\sigma\!\mu_v\in \prod_{v \in S_\infty}X^+(T_n)$ defined by ${}^\sigma\!\mu_v = \mu_{\sigma^{-1}\circ v}$ for $v \in S_\infty$. We have the $\sigma$-linear isomorphism $M_{\mu,\C} \rightarrow M_{{}^\sigma\!\mu,\C}$ defined by
\begin{align}\label{E:sigma linear}
\bigotimes_{v \in S_\infty}z_v\cdot{\bf v}_{v}= \bigotimes_{v \in S_\infty} \sigma(z_v) \cdot {\bf v}_{\sigma^{-1}\circ v}
\end{align}
for $z_v \in \C$, ${\bf v}_{v} \in M_{\mu_v}$.

\subsubsection{Irreducible algebraic representations of $\GL_2$}
Let $\mu = (\mu_1,\mu_2) \in X^+(T_2)$.
We fix a model of the representation $(\rho_\mu,M_\mu)$ as follows:
Let $M_\mu$ be the $\Q$-vector space consisting of homogeneous polynomials over $\Q$ of degree $\mu_1-\mu_2$ in variables $x$ and $y$. For $g \in \GL_2(\Q)$, we define $\rho_\mu(g) \in {\rm End}_\Q(M_\mu)$ by
\[
\rho_\mu(g)\cdot P(x,y) = \det(g)^{\mu_2}\cdot P((x,y)g)
\]
for $P \in M_\mu$.

\subsubsection{Irreducible algebraic representations of $\GL_3$}

Let $\mu = (\mu_1,\mu_2,\mu_3) \in X^+(T_3)$. We fix a model of the representation $(\rho_\mu,M_\mu)$ as follows:
Let $M_\mu$ be the $\Q$-vector space generated by homogeneous polynomials of degree $\mu_1+\mu_2-2\mu_3$ in variable $x = (x_{ij})_{1 \leq i \leq 2,\, 1 \leq j \leq 3}$ of the form
\[
\prod_{i=1}^3 x_{2i}^{n_i} \prod_{1 \leq j < j' \leq 3}\det \bp x_{1j} & x_{1j'} \\ x_{2j} & x_{2j'}\ep^{n_{j,j'}}
\]
with $n_1+n_2+n_3 = \mu_1-\mu_2$ and $n_{12}+n_{23}+n_{13} = \mu_2-\mu_3$.
For $g \in \GL_3(\Q)$, we define $\rho_\mu(g) \in {\rm End}_\Q(M_\mu)$ by
\[
\rho_\mu(g)\cdot P(x) = \det(g)^{\mu_3}\cdot P(xg)
\]
for $P \in M_\mu$.
Let $P_\mu^+$ be the highest weight vector of $M_\mu$ given by
\[
P_{\mu}^+=x_{21}^{\mu_1-\mu_2}\det \bp x_{11} & x_{12} \\ x_{21} & x_{22} \ep^{\mu_2-\mu_3}.
\]

\subsubsection{Irreducible algebraic representations of ${\rm SO}(3)$}

Let $\ell \in \Z_{\geq 0}$. Let $(\tau_\ell,V_\ell)$ be the irreducible algebraic representation of ${\rm SO}(3)$ of dimension $2\ell+1$. We fix a model of the representation as follows:
Let $V_\ell$ be the quotient of the space of homogeneous polynomials over $\C$ of degree $\ell$ in variables $x_1,x_2,x_3$ by the subspace generated by $x_1^2+x_2^2+x_3^2$.
Denote by $\overline{x}_1,\overline{x}_2,\overline{x}_3$ the image of $x_1,x_2,x_3$ in $V_\ell$.
For $g \in {\rm SO}(3)$, we define $\tau_\ell(g) \in {\rm End}_\C(V_\ell)$ by
\[
\tau_\ell(g)\cdot P(\overline{x}_1,\overline{x}_2,\overline{x}_3) = P((\overline{x}_1,\overline{x}_2,\overline{x}_3)g)
\]
for $P \in V_\ell$.
For $(j_1,j_2,j_3) \in \Z_{\geq 0}^3$ with $j_1+j_2+j_3 = \ell$, let
\[
{\bf v}_{(\ell;\,(j_1,j_2,j_3))} = \overline{x}_1^{j_1}\overline{x}_2^{j_2}\overline{x}_3^{j_3}.
\]
Let $\{{\bf v}_{(\ell;\,i)} \, \vert \, -\ell \leq i \leq \ell\}$ be the basis of $V_\ell$ defined by
\[
{\bf v}_{(\ell;\,i)} = ({\rm sgn}(i)\overline{x}_1 + \sqrt{-1}\,\overline{x}_2)^{|i|}\overline{x}_3^{\ell-|i|}
\]
for $-\ell \leq i \leq \ell$.
We have the following relations 
\begin{align}\label{E:Lie algebra action}
\begin{split}
E_{12}\cdot {\bf v}_{(\ell;\,i)} &= \sqrt{-1}\,i\cdot {\bf v}_{(\ell;\,i)},\\
E_\pm \cdot {\bf v}_{(\ell;\,i)} &= (\pm \ell -i)\cdot {\bf v}_{(\ell;\,i\pm 1)}
\end{split}
\end{align}
for $-\ell \leq i \leq \ell$, where $\{E_+,E_{12},E_-\}$ is the basis of $\frak{so}(3)_\C$ defined by
\[
E_{12} = \bp 0&1&0 \\ -1&0&0 \\ 0&0&0\ep,\quad E_\pm = \pm\bp 0&0&0 \\ 0&0&1 \\ 0&-1&0\ep \otimes \sqrt{-1} + \bp 0&0&1\\0&0&0\\-1&0&0\ep.
\]
For example, the adjoint representation of ${\rm SO}(3)$ on $\frak{g}_{3,\C} / \frak{k}_{3,\C}$ is isomorphic to $V_{2}$.
Moreover, the basis $\{X_2,X_1,X_0,X_{-1},X_{-2}\}$ in (\ref{E:Lie algebra basis}) satisfies the analogue relations (\ref{E:Lie algebra action}), that is, 
\begin{align}\label{E:Lie algebra action 2}
\begin{split}
{\rm ad}(E_{12})\cdot X_i &= \sqrt{-1}\,i\cdot X_i,\\
{\rm ad}(E_\pm) \cdot X_i &= (\pm 2 -i)\cdot X_{i\pm 1}
\end{split}
\end{align}
for $-2 \leq i \leq 2$.

\subsection{Rational structures via the cuspidal cohomology}\label{SS:cohomology}
Let $\mu = \prod_{v \in S_\infty} \mu_v \in \prod_{v \in S_\infty} X^+(T_3)$. 
Let $K_f$ be a neat open compact subgroup of $\GL_3(\A_{\F,f})$ and $\mathcal{S}_{3,K_f}$ be the locally symmetric space defined by
\[
\mathcal{S}_{3,K_f} = \GL_3(\F)\backslash \GL_3(\A_\F) / K_3K_f.
\]
The irreducible algebraic representation $M_{\mu,\C}^\vee$ defines a locally constant sheaf
$\mathcal{M}_{\mu,\C}^\vee$ of $\C$-vector spaces on $\mathcal{S}_{3,K_f}$. For $q \in \Z_{\geq 0}$, we denote by
\[
H^q_B(\mathcal{S}_{3,K_f},{M}_{\mu,\C}^\vee),\quad H^q(\mathcal{S}_{3,K_f},\mathcal{M}_{\mu,\C}^\vee)
\]
the $q$-th singular cohomology of $\mathcal{S}_{3,K_f}$ with coefficients in ${M}_{\mu,\C}^\vee$ and the $q$-th sheaf cohomology of $\mathcal{M}^\vee_{\mu,\C}$, respectively.
The two cohomology groups are canonically isomorphic under the de Rham isomorphism.
For $\sigma \in {\rm Aut}(\C)$, the canonical $\sigma$-linear isomorphism $M_{\mu,\C} \rightarrow M_{{}^\sigma\!\mu,\C}$ in (\ref{E:sigma linear}) induces a $\sigma$-linear isomorphism from $H^q_B(\mathcal{S}_{3,K_f},{M}_{\mu,\C}^\vee)$ to $H^q_B(\mathcal{S}_{3,K_f},{M}_{{}^\sigma\!\mu,\C}^\vee)$. This in term induces a $\sigma$-linear isomorphism
\[
T_{\sigma,K_f} : H^q(\mathcal{S}_{3,K_f},\mathcal{M}_{\mu,\C}^\vee) \longrightarrow H^q(\mathcal{S}_{3,K_f},\mathcal{M}_{{}^\sigma\!\mu,\C}^\vee)
\]
via the de Rham isomorphism.
Passing to the limit and define
\[
H^q(\mathcal{S}_3,\mathcal{M}_{\mu,\C}^\vee) = \varinjlim_{K_f}H^q(\mathcal{S}_{3,K_f},\mathcal{M}_{\mu,\C}^\vee).
\]
Then $H^q(\mathcal{S}_3,\mathcal{M}_{\mu,\C}^\vee)$ is canonically isomorphic to the relative Lie algebra cohomology
\[
H^q(\frak{g}_{3,\infty},K_{3,\infty};C^{\infty}(\GL_3(\F)\backslash \GL_3(\A_\F)) \otimes M_{\mu,\C}^\vee).
\]
We have the $\GL_3(\A_{\F,f})$-module structure on $H^q(\mathcal{S}_3,\mathcal{M}_{\mu,\C}^\vee)$ induced by the right translation of $\GL_3(\A_{\F,f})$ on $C^{\infty}(\GL_3(\F)\backslash \GL_3(\A_\F))$.
For $\sigma \in {\rm Aut}(\C)$, we have the $\sigma$-linear isomorphism 
\[
T_\sigma = \varinjlim_{K_f}T_{\sigma,K_f} : H^q(\mathcal{S}_3,\mathcal{M}_{\mu,\C}^\vee) \longrightarrow H^q(\mathcal{S}_3,\mathcal{M}_{{}^\sigma\!\mu,\C}^\vee)
\]
which commutes with the $\GL_3(\A_{\F,f})$-action.
The cuspidal cohomology is define to be the relative Lie algebra cohomology
\[
H_{\rm cusp}^q(\mathcal{S}_3,\mathcal{M}_{\mu,\C}^\vee) = H^q(\frak{g}_{3,\infty},K_{3,\infty};\mathcal{A}_{\rm cusp}(\GL_3(\F)\backslash \GL_3(\A_\F)) \otimes M_{\mu,\C}^\vee),
\]
where $\mathcal{A}_{\rm cusp}(\GL_3(\F)\backslash \GL_3(\A_\F))$ is the space of cusp forms on $\GL_3(\A_\F)$.
As explained in \cite[p.\,126]{Clozel1990}, the natural inclusion 
\[
\mathcal{A}_{\rm cusp}(\GL_3(\F)\backslash \GL_3(\A_\F)) \subset C^\infty(\GL_3(\F)\backslash \GL_3(\A_\F))
\]
induces a $\GL_3(\A_{\F,f})$-equivariant injective homomorphism
\[
H_{\rm cusp}^q(\mathcal{S}_3,\mathcal{M}_{\mu,\C}^\vee) \longrightarrow H^q(\mathcal{S}_3,\mathcal{M}_{\mu,\C}^\vee).
\]
We listed in the following theorem some general results on the cuspidal cohomology for $\GL_3$. Similar assertions hold for $\GL_n$ and the proofs are contained in \cite[\S\,3.5 and Lemme 4.9]{Clozel1990}.

\begin{thm}\label{T:Clozel}
Let $\mu = \prod_{v \in S_\infty} \mu_v \in \prod_{v \in S_\infty} X^+(T_3)$ and $q \in \Z_{\geq 0}$.
\begin{itemize}
\item[(1)] For $\sigma\in{\rm Aut}(\C)$, we have $T_\sigma\left(H_{\rm cusp}^q(\mathcal{S}_3,\mathcal{M}_{\mu,\C}^\vee)\right) = H_{\rm cusp}^q(\mathcal{S}_3,\mathcal{M}_{{}^\sigma\!\mu,\C}^\vee)$.
\item[(2)] $H_{\rm cusp}^q(\mathcal{S}_3,\mathcal{M}_{\mu,\C}^\vee)=0$ unless $\mu$ is pure and $2d \leq q \leq 3d$.
\item[(3)] Suppose $\mu$ is pure and write
\[
\mu_v = (\tfrac{\ell_v-3+{\sf w}}{2},\tfrac{{\sf w}}{2},\tfrac{-\ell_v+3+{\sf w}}{2})
\]
for some odd integer $\ell_v \in \Z_{\geq 3}$ and even integer ${\sf w}$ for $v \in S_\infty$. 
Then
\[
H_{\rm cusp}^q(\mathcal{S}_3,\mathcal{M}_{\mu,\C}^\vee) = \bigoplus_\itSigma H^q(\frak{g}_{3,\infty},K_{3,\infty};\itSigma \otimes M_{\mu,\C}^\vee),
\]
where $\itSigma$ runs through irreducible cuspidal automorphic representations of $\GL_3(\A_\F)$ such that
\[
\itSigma_v = {\rm Ind}^{\GL_3(\R)}_{P_{2,1}(\R)}(D_{\ell_v} \boxtimes \epsilon_v) \otimes |\mbox{ }|^{{\sf w}/2}
\]
for some quadratic character $\epsilon_v$ of $\R^\times$, $v \in S_\infty$.
Here $P_{2,1}$ is the standard parabolic subgroup of $\GL_3$ with Levi part $\GL_2 \times \GL_1$ and $D_\ell$ is the discrete series representation of $\GL_2(\R)$ of weight $\ell$.
\item[(4)] Let $\itSigma$ be an irreducible cuspidal automorphic representation which contributes to the cuspidal cohomology of $\GL_3(\A_\F)$ with coefficients in $M^\vee_{\mu,\C}$. 
Then
\[
H^q(\frak{g}_{3,\infty},K_{3,\infty};\itSigma \otimes M_{\mu,\C}^\vee) = \left(\itSigma \otimes \bigwedge^q (\frak{g}_{3,\infty,\C} / \frak{k}_{3,\infty,\C})^* \otimes M_{\mu,\C}^\vee\right)^{K_{3,\infty}} \simeq \itSigma_f
\]
for $q=2d,3d$.
\end{itemize}
\end{thm}

We call $2d$ and $3d$ the bottom and top degree cuspidal degree, respectively.
Let $\itSigma = \bigotimes_v \itSigma_v$ be an irreducible cuspidal automorphic representation which contributes to the cuspidal cohomology of $\GL_3(\A_\F)$ with coefficients in $M^\vee_{\mu,\C}$. Let $\itSigma_f  = \bigotimes_{v \nmid \infty} \itSigma_v$ be the finite part of $\itSigma$. For $\sigma \in {\rm Aut}(\C)$, let ${}^\sigma\!\itSigma$ be the irreducible admissible representation of $\GL_3(\A_\F)$ defined by 
\[
{}^\sigma\!\itSigma = \bigotimes_{v \in S_\infty} {}^\sigma\!\itSigma_v \otimes {}^\sigma\!\itSigma_f,
\]
where ${}^\sigma\!\itSigma_v = \itSigma_{\sigma^{-1}\circ v}$ for $v \in S_\infty$.
By Theorem \ref{T:Clozel} and the strong multiplicity one theorem for $\GL_3$, ${}^\sigma\!\itSigma$ is automorphic cuspidal and contributes to the cuspidal cohomology of $\GL_3(\A_\F)$ with coefficients in $M^\vee_{{}^\sigma\!\mu,\C}$.
Moreover, ${}^\sigma\!\itSigma_f$ appears in the cuspidal cohomology $H_{\rm cusp}^q(\mathcal{S}_3,\mathcal{M}_{\mu,\C}^\vee)$ with multiplicity one for $q=2d,3d$, and the ${}^\sigma\!\itSigma_f$-isotypic component is equal to $H^q(\frak{g}_{3,\infty},K_{3,\infty};{}^\sigma\!\itSigma \otimes M_{{}^\sigma\!\mu,\C}^\vee)$. We have the $\sigma$-linear $\GL_3(\A_{\F,f})$-equivariant isomorphism 
\begin{align}\label{E:sigma linear map 1}
T_\sigma : H^q(\frak{g}_{3,\infty},K_{3,\infty}; \itSigma \otimes M_{\mu,\C}^\vee) \longrightarrow H^q(\frak{g}_{3,\infty},K_{3,\infty};{}^\sigma\!\itSigma \otimes M_{{}^\sigma\!\mu,\C}^\vee).
\end{align}
The rationality field $\Q(\itSigma)$ of $\itSigma$ is define to be the fixed field of
$\left\{\sigma \in {\rm Aut}(\C) \, \vert \, {}^\sigma\!\itSigma = \itSigma \right\}$ and is a number field.
We then have the $\Q(\itSigma)$-rational structure on the cuspidal cohomology $H^q(\frak{g}_{3,\infty},K_{3,\infty}; \itSigma \otimes M_{\mu,\C}^\vee)$ by taking the ${\rm Aut}(\C/\Q(\itSigma))$-invariants
\begin{align*}
&H^q(\frak{g}_{3,\infty},K_{3,\infty}; \itSigma \otimes M_{\mu,\C}^\vee)^{{\rm Aut}(\C/\Q(\itSigma))}\\
& = \left.\left\{\varphi \in H^q(\frak{g}_{3,\infty},K_{3,\infty}; \itSigma \otimes M_{\mu,\C}^\vee) \,\right\vert\, T_\sigma \varphi=\varphi  \mbox{ for all }\sigma \in {\rm Aut}(\C/\Q(\itSigma))\right\}.
\end{align*}

\subsection{Rational structures via the Whittaker model}\label{SS:Whittaker}

Let $\itSigma = \bigotimes_v \itSigma_v$ be a cohomological irreducible cuspidal automorphic representation of $\GL_3(\A_\F)$ with central character $\omega_\itSigma$. 
For each place $v$ of $\F$, let $\mathcal{W}(\itSigma_v,\psi_{3,\F_v})$ be the space of Whittaker functions of $\itSigma_v$ with respect to $\psi_{3,\F_v}$. Recall that $\mathcal{W}(\itSigma_v,\psi_{3,\F_v})$ is contained in the space of locally constant (resp.\,smooth and moderate growth) functions $W : \GL_3(\F_v) \rightarrow \C$ which satisfy
\[
W(ug) = \psi_{3,\F_v}(u)W(g)
\]
for all $u \in N_3(\F_v)$ and $g \in \GL_3(\F_v)$ when $v$ is finite (resp.\,$v$ is archimedean).
Let 
\[
\mathcal{W}(\itSigma_v,\psi_{3,\F_\infty}) = \bigotimes_{v \in S_\infty} \mathcal{W}(\itSigma_v,\psi_{n,\R}).
\]
Let $\mathcal{W}(\itSigma_f,\psi_{3,\F}^{(\infty)})$ be the space of Whittaker function of $\itSigma_f = \bigotimes_{v \nmid \infty} \itSigma_v$ with respect to $\psi_{3,\F}^{(\infty)} = \bigotimes_{v \nmid \infty} \psi_{3,\F_v}$. 
We have the $\GL_3(\A_{\F,f})$-equivariant isomorphism
\[
\bigotimes_{v \nmid \infty}\mathcal{W}(\itSigma_v,\psi_{3,\F_v}) \longrightarrow \mathcal{W}(\itSigma_f,\psi_{3,\F}^{(\infty)}),\quad \bigotimes_{v \nmid \infty} W_v \longmapsto \prod_{v \nmid \infty} W_v.
\]
Here we take the restricted tensor product on the left-hand side with respect to the normalized unramified Whittaker function $W_v^\circ \in \mathcal{W}(\itSigma_v,\psi_{3,\F_v})$, that is, the right $\GL_3(\frak{o}_{\F_v})$-invariant Whittaker function normalized so that $W_v^\circ(1)=1$ for all but finitely many finite places $v$ such that $\itSigma_v$ is unramified.
For $\sigma \in {\rm Aut}(\C)$, let 
\begin{align}\label{E:sigma linear map 2}
t_\sigma : \mathcal{W}(\itSigma_f,\psi_{3,\F}^{(\infty)}) \longrightarrow \mathcal{W}({}^\sigma\!\itSigma_f,\psi_{3,\F}^{(\infty)})
\end{align}
be the $\sigma$-linear $\GL_3(\A_{\F,f})$-equivariant isomorphism defined by
\[
t_\sigma W(g) = \sigma ( W ({\rm diag}(u_\sigma^{-2},u_\sigma^{-1},1)g))
\]
for $g \in \GL_3(\A_{\F,f})$. Here $u_\sigma \in  \prod_p \Z_p^\times \subset \A_{\F,f}^\times$ is the unique element depending on $\sigma$ such that $\sigma(\psi_\Q(x)) = \psi_\Q(u_\sigma x)$ for all $x \in \A_{\Q,f}$.
We then have the $\Q(\itSigma)$-rational structure on $\mathcal{W}(\itSigma_f,\psi_{3,\F}^{(\infty)})$ by taking the ${\rm Aut}(\C/\Q(\itSigma))$-invariants 
\[
\mathcal{W}(\itSigma_f,\psi_{3,\F}^{(\infty)})^{{\rm Aut}(\C/\Q(\itSigma))} = \left.\left\{W \in \mathcal{W}(\itSigma_f,\psi_{3,\F}^{(\infty)}) \,\right\vert\, t_\sigma W=W  \mbox{ for all }\sigma \in {\rm Aut}(\C/\Q(\itSigma))\right\}.
\]

Let $\mu = \prod_{v \in S_\infty} \mu_v\in \prod_{v \in S_\infty}X^+(T_3)$ such that $\itSigma$ contributes to the cuspidal cohomology of $\GL_3(\A_\F)$ with coefficients in $M_{\mu,\C}^\vee$ with central character $\omega_\itSigma$.
We have $|\omega_\itSigma| = |\mbox{ }|_{\A_\F}^{{\sf w}}$ for some even integer ${\sf w}$. 
For each $v \in S_\infty$, we have
\[
\mu_v = (\tfrac{\ell_v-3+{\sf w}}{2},\tfrac{{\sf w}}{2},\tfrac{-\ell_v+3+{\sf w}}{2})
\]
for some odd integer $\ell_v \in \Z_{\geq 3}$ and 
\[
\itSigma_v = {\rm Ind}^{\GL_3(\R)}_{P_{2,1}(\R)}(D_{\ell_v} \boxtimes \epsilon_v) \otimes |\mbox{ }|^{{\sf w}/2}
\]
for some quadratic character $\epsilon_v$ of $\R^\times$. 
Note that (cf.\,\cite[Corollary 3.6]{Miyazaki2009})
\[
\itSigma_v \vert_{{\rm SO}(3)} = \bigoplus_{i=0}^\infty
(i+1)\cdot (V_{\ell_v+2i}\oplus V_{\ell_v+2i+1}).
\] 
In particular, $V_{\ell_v}$ appears in $\itSigma_v\vert_{{\rm SO}(3)}$ with multiplciity one. We call $V_{\ell_v}$ the minimal ${\rm SO}(3)$-type of $\itSigma_v$.
In the following lemma, we write down generators of the one-dimensional vector spaces $H^{2d}(\frak{g}_{3,\infty},K_{3,\infty};\itSigma_\infty \otimes M_{\mu,\C}^\vee)$ and  $H^{3d}(\frak{g}_{3,\infty},K_{3,\infty};\itSigma_\infty \otimes M_{\mu,\C}^\vee)$.
\begin{lemma}\label{L:cohomology generator}
Let $v \in S_\infty$.
The vectors
\[
\sum_{i=-\ell_v}^{\ell_v}\frac{1}{(\ell_v+i)!}\cdot {\bf v}_{(\ell_v;\,i)}\otimes E_-^{\ell_v+i}\cdot \left(X_{-1}^*\wedge X_{-2}^* \otimes \rho_{\mu_v^\vee}\left( \bp 1 & 0 & 1  \\ \sqrt{-1} & 0 & -\sqrt{-1} \\ 0 & 1 & 0\ep\right)P_{\mu_v^\vee}^+\right)
\]
and
\[
\sum_{i=-\ell_v}^{\ell_v}\frac{1}{(\ell_v+i)!}\cdot {\bf v}_{(\ell_v;\,i)}\otimes E_-^{\ell_v+i}\cdot \left(X_0^*\wedge X_{-1}^*\wedge X_{-2}^* \otimes \rho_{\mu_v^\vee}\left( \bp 1 & 0 & 1  \\ \sqrt{-1} & 0 & -\sqrt{-1} \\ 0 & 1 & 0\ep\right)P_{\mu_v^\vee}^+\right)
\]
are generators of $\left(\itSigma_v\otimes \bigwedge^2 (\frak{g}_{3,\C} / \frak{k}_{3,\C})^* \otimes M_{\mu_v,\C}^\vee\right)^{K_3}$ and $\left(\itSigma_v\otimes \bigwedge^3 (\frak{g}_{3,\C} / \frak{k}_{3,\C})^* \otimes M_{\mu_v,\C}^\vee\right)^{K_3}$, respectively.
Moreover, we have
\begin{align}\label{E:relation}
\begin{split}
&\frac{1}{(\ell_v+i)!}\cdot E_-^{\ell_v+i}\cdot \left(X_{-1}^*\wedge X_{-2}^* \otimes \rho_{\mu_v^\vee}\left( \bp 1 & 0 & 1  \\ \sqrt{-1} & 0 & -\sqrt{-1} \\ 0 & 1 & 0\ep\right)P_{\mu_v^\vee}^+\right)\\
& = \frac{(-1)^{i+{{\sf w}(\itSigma)}/2}}{(\ell_v-i)!}\cdot E_+^{\ell_v-i}\cdot \left(X_1^*\wedge X_2^* \otimes \rho_{\mu_v^\vee}\left( \bp -1 & 0 & -1  \\ \sqrt{-1} & 0 & -\sqrt{-1} \\ 0 & 1 & 0\ep\right)P_{\mu_v^\vee}^+\right)
\end{split}
\end{align}
for $-\ell_v \leq i \leq \ell_v$. Similar formula holds if we replace $X_{-1}^*\wedge X_{-2}^*$ and $X_1^*\wedge X_2^*$ by $X_0^*\wedge X_{-1}^*\wedge X_{-2}^*$ and $X_0^*\wedge X_1^*\wedge X_2^*$, respectively.
\end{lemma}

\begin{proof}
We drop the subscript $v$ for brevity.
Firstly we show that
\begin{align}\label{E:fixed vector}
(V_\ell \otimes V_\ell)^{{\rm SO}(3)} = \sum_{i=-\ell}^{\ell}\frac{(-1)^i}{(\ell-i)!(\ell+i)!}\cdot {\bf v}_{(\ell;\,i)}\otimes {\bf v}_{(\ell;\,-i)}\cdot \C.
\end{align}
Let $\sum_{-\ell \leq i,j \leq \ell}c_{ij}\cdot {\bf v}_{(\ell;\,i)}\otimes {\bf v}_{(\ell;\,j)}$ be a non-zero vector in $(V_\ell \otimes V_\ell)^{{\rm SO}(3)}$. Since
\[
E_{12}\cdot \sum_{-\ell \leq i,j \leq \ell}c_{ij}\cdot {\bf v}_{(\ell;\,i)}\otimes {\bf v}_{(\ell;\,j)} = \sqrt{-1}\sum_{-\ell \leq i,j \leq \ell}(i+j)c_{ij}\cdot {\bf v}_{(\ell;\,i)}\otimes {\bf v}_{(\ell;\,j)},
\]
we see that $c_{ij}=0$ unless $i+j=0$. Thus 
\[
\sum_{-\ell \leq i,j \leq \ell}c_{ij}\cdot {\bf v}_{(\ell;\,i)}\otimes {\bf v}_{(\ell;\,j)} = \sum_{i=-\ell}^{\ell} c_{i}\cdot {\bf v}_{(\ell;\,i)}\otimes {\bf v}_{(\ell;\,-i)}
\]
for some $c_i$. By (\ref{E:Lie algebra action}), we have
\[
E_+\cdot \sum_{i=-\ell}^{\ell} c_{i}\cdot {\bf v}_{(\ell;\,i)}\otimes {\bf v}_{(\ell;\,-i)} = \sum_{i=-\ell}^{\ell-1} [(\ell-i)c_i+(\ell+i+1)c_{i+1}]\cdot{\bf v}_{(\ell;\,i)}\otimes {\bf v}_{(\ell;\,-i)}.
\]
We thus obtain the recursive relation 
$
(\ell-i)c_i+(\ell+i+1)c_{i+1}=0
$
for $-\ell \leq i \leq \ell-1$, from which (\ref{E:fixed vector}) follows immediately. By (\ref{E:Lie algebra action}) again,
\[
E_-^{\ell+i}\cdot {\bf v}_{(\ell;\,\ell)} = (-1)^{\ell+i}\frac{(2\ell)!}{(\ell-i)!}\cdot {\bf v}_{(\ell;\,-i)}.
\]
Thus we can rewrite (\ref{E:fixed vector}) as
\begin{align}\label{E:fixed vector 2}
(V_\ell \otimes V_\ell)^{{\rm SO}(3)} = \sum_{i=-\ell}^{\ell}\frac{1}{(\ell+i)!}\cdot {\bf v}_{(\ell;\,i)}\otimes E_-^{\ell+i}\cdot {\bf v}_{(\ell;\,\ell)}\cdot \C.
\end{align}
Now we regard ${\rm SO}(2)$ as a Cartan subgroup of ${\rm SO}(3)$ via the embedding $g \mapsto \bp g & 0 \\ 0 & 1 \ep$.
Note that 
\[
\bigwedge^2 (\frak{g}_{3,\C} / \frak{k}_{3,\C})^* \simeq \bigwedge^3 (\frak{g}_{3,\C} / \frak{k}_{3,\C})^* \simeq  V_1 \oplus V_3
\]
and $V_{\ell-3}$ appears in $M_{\mu,\C}^\vee \vert_{{\rm SO}(3)}$ with multiplicity one.
Moreover, it is easy to see that the highest weight vectors (with respect to ${\rm SO}(2)$) of $\bigwedge^2 (\frak{g}_{3,\C} / \frak{k}_{3,\C})^*$, $\bigwedge^3 (\frak{g}_{3,\C} / \frak{k}_{3,\C})^*$, and $M_{\mu,\C}^\vee \vert_{{\rm SO}(3)}$ are
\[
X_{-1}^* \wedge X_{-2}^*,\quad X_0^*\wedge X_{-1}^* \wedge X_{-2}^*,\quad \rho_{\mu^\vee}\left( \bp 1 & 0 & 1  \\ \sqrt{-1} & 0 & -\sqrt{-1} \\ 0 & 1 & 0\ep\right)P_{\mu^\vee}^+
\]
respectively.
Therefore, we have ${\rm SO}(3)$-equivariant homomorphisms 
\begin{align*}
V_\ell &\longrightarrow \bigwedge^2 (\frak{g}_{3,\C} / \frak{k}_{3,\C})^* \otimes M_{\mu,\C}^\vee, \quad {\bf v}_{(\ell;\,\ell)} \longmapsto X_{-1}^* \wedge X_{-2}^* \otimes \rho_{\mu^\vee}\left( \bp 1 & 0 & 1  \\ \sqrt{-1} & 0 & -\sqrt{-1} \\ 0 & 1 & 0\ep\right)P_{\mu^\vee}^+,\\
V_\ell &\longrightarrow \bigwedge^3 (\frak{g}_{3,\C} / \frak{k}_{3,\C})^* \otimes M_{\mu,\C}^\vee, \quad {\bf v}_{(\ell;\,\ell)} \longmapsto X_0^*\wedge X_{-1}^* \wedge X_{-2}^* \otimes \rho_{\mu^\vee}\left( \bp 1 & 0 & 1  \\ \sqrt{-1} & 0 & -\sqrt{-1} \\ 0 & 1 & 0\ep\right)P_{\mu^\vee}^+.
\end{align*}
The first assertion then follows from (\ref{E:fixed vector 2}). To prove the second assertion, note that ${\rm O}(3)/{\rm SO}(3)$ acts on $H^2(\frak{g}_3,K_3;\itSigma_\infty \otimes M_{\mu,\C}^\vee)$ and  $H^3(\frak{g}_3,K_3;\itSigma_\infty \otimes M_{\mu,\C}^\vee)$ by a sign. 
By taking the representative ${\rm diag}(-1,-1,-1)$, we see that the sign is equal to $(-1)^{{\sf w}(\itSigma)/2}\omega_{\itSigma_\infty}(-1)$. On the other hand, taking the representative ${\rm diag}(-1,1,1)$, we have
\[
{\rm diag}(-1,1,1)\cdot {\bf v}_{(\ell;\,i)} = -\omega_{\itSigma_\infty}(-1)\cdot {\bf v}_{(\ell;\,-i)}
\]
in $\itSigma_\infty$ and 
\begin{align*}
&{\rm diag}(-1,1,1)\cdot E_-^{\ell+i}\cdot \left(X_{-1}^*\wedge X_{-2}^* \otimes \rho_{\mu^\vee}\left( \bp 1 & 0 & 1  \\ \sqrt{-1} & 0 & -\sqrt{-1} \\ 0 & 1 & 0\ep\right)P_{\mu^\vee}^+\right)\\
& = (-1)^{1+i} E_+^{\ell-i}\cdot \left(X_1^*\wedge X_2^* \otimes \rho_{\mu^\vee}\left( \bp -1 & 0 & -1  \\ \sqrt{-1} & 0 & -\sqrt{-1} \\ 0 & 1 & 0\ep\right)P_{\mu^\vee}^+\right).
\end{align*}
Similar formula holds if we replace $X_{-1}^*\wedge X_{-2}^*$ and $X_1^*\wedge X_2^*$ by $X_0^*\wedge X_{-1}^*\wedge X_{-2}^*$ and $X_0^*\wedge X_1^*\wedge X_2^*$, respectively.
Thus (\ref{E:relation}) follows immediately. This completes the proof.
\end{proof}

Let $v \in S_\infty$.
The following theorem of Miyazaki \cite[Theorem 5.9]{Miyazaki2009} (see also \cite[Theorem 3.1]{HIM2012}) is crucial to our normalization of generators of $H^{2}(\frak{g}_{3},K_{3};\itSigma_v \otimes M_{\mu_v,\C}^\vee)$ and  $H^{3}(\frak{g}_{3},K_{3};\itSigma_v \otimes M_{\mu_v,\C}^\vee)$.
\begin{thm}[Miyazaki]\label{T:Miyazaki}
There exists a unique ${\rm SO}(3)$-equivariant homomorphism 
\[
V_{\ell_v} \longrightarrow \mathcal{W}(\itSigma_v,\psi_{3,\R}),\quad {\bf v}_{(\ell_v;\, (j_1,j_2,j_3))} \longmapsto W_{(\ell_v, {\sf w}, \epsilon_v;\, (j_1,j_2,j_3))}
\]
such that $W_{(\ell_v, {\sf w}, \epsilon_v;\, (j_1,j_2,j_3))}(g\cdot {\rm diag}(-1,-1,-1)) = (-1)^{1+{\sf w}/2}\epsilon_v(-1)\cdot W_{(\ell_v, {\sf w}, \epsilon_v;\,(j_1,j_2,j_3))}(g)$ and 
\begin{align*}
W_{(\ell_v, {\sf w}, \epsilon_v;\,(j_1,j_2,j_3))}({\rm diag}(a_1a_2a_3,a_2a_3,a_3)) &= (\sqrt{-1})^{j_1-j_3} (a_1a_2^2a_3^3)^{{\sf w}/2}\int_{L_1}\frac{ds_1}{2\pi\sqrt{-1}}\int_{L_2}\frac{ds_2}{2\pi\sqrt{-1}}\,a_1^{-s_1+1}a_2^{-s_2+1}\\
&\quad\quad\quad\quad\quad\times \frac{\Gamma_\C(s_1+\tfrac{\ell_v-1}{2})\Gamma_\R(s_1+j_1)\Gamma_\C(s_2+\tfrac{\ell_v-1}{2})\Gamma_\R(s_2+j_3)}{\Gamma_\R(s_1+s_2+j_1+j_3)}
\end{align*}
for all $a_1,a_2, a_3>0$ and $(j_1,j_2,j_3) \in \Z_{\geq 0}^3$ with $j_1+j_2+j_3=\ell_v$,
where $L_1$ and $L_2$ are any vertical paths from south to north which keep the poles of $\Gamma_\C(s_1+\tfrac{\ell_v-1}{2})\Gamma_\R(s_1+j_1)$ and $\Gamma_\C(s_2+\tfrac{\ell_v-1}{2})\Gamma_\R(s_2+j_3)$, respectively, on its left.
\end{thm}

For $-\ell_v \leq i \leq \ell_v$, let 
\[
W_{(\ell_v, {\sf w},\epsilon_v;\,i)} \in \mathcal{W}(\itSigma_v,\psi_{3,\R})
\]
be the image of ${\bf v}_{(\ell_v;\, i)}$ under the homomorphism in Theorem \ref{T:Miyazaki}. 
Define the classes 
\[
[\itSigma_v]_b  \in H^2(\frak{g}_3,K_3;\mathcal{W}(\itSigma_v,\psi_{3,\R}) \otimes M_{\mu_v,\C}^\vee),\quad [\itSigma_v]_t  \in H^3(\frak{g}_3,K_3;\mathcal{W}(\itSigma_v,\psi_{3,\R}) \otimes M_{\mu_v,\C}^\vee)
\]
by 
\begin{align*}
\begin{split}
[\itSigma_v]_b
&= \sum_{i=-\ell_v}^{\ell_v}\frac{(\sqrt{-1})^{{\sf w}/2}}{(\ell_v+i)!}\cdot W_{(\ell_v,{\sf w}, \epsilon_v;\,i)}\otimes E_-^{\ell_v+i}\cdot \left(X_{-1}^*\wedge X_{-2}^* \otimes \rho_{\mu_v^\vee}\left( \bp 1 & 0 & 1  \\ \sqrt{-1} & 0 & -\sqrt{-1} \\ 0 & 1 & 0\ep\right)P_{\mu_v^\vee}^+\right),\\
[\itSigma_v]_t
&= \sum_{i=-\ell_v}^{\ell_v}\frac{(\sqrt{-1})^{1+{\sf w}/2}}{(\ell_v+i)!}\cdot W_{(\ell_v,{\sf w}, \epsilon_v;\,i)}\otimes E_-^{\ell_v+i}\cdot \left(X_0^*\wedge X_{-1}^*\wedge X_{-2}^* \otimes \rho_{\mu_v^\vee}\left( \bp 1 & 0 & 1  \\ \sqrt{-1} & 0 & -\sqrt{-1} \\ 0 & 1 & 0\ep\right)P_{\mu_v^\vee}^+\right).
\end{split}
\end{align*}
Note that the classes are well-define by Lemma \ref{L:cohomology generator} and the subscripts $b$ and $t$ refer to classes in the bottom and top degree cuspidal cohomology.
Now we define
\begin{align}\label{E:GL_3 class}
\begin{split}
[\itSigma_\infty]_b = \bigotimes_{v \in S_\infty}[\itSigma_v]_b &\in H^{2d}\left(\frak{g}_{3,\infty},K_{3,\infty};\mathcal{W}(\itSigma_\infty,\psi_{3,\F_\infty}) \otimes M_{\mu,\C}^\vee\right)\\
& = \bigotimes_{v \in S_\infty}H^2(\frak{g}_3,K_3;\mathcal{W}(\itSigma_v,\psi_{2,\R}) \otimes M_{\mu_v,\C}^\vee),\\
[\itSigma_\infty]_t = \bigotimes_{v \in S_\infty}[\itSigma_v]_t &\in H^{3d}\left(\frak{g}_{3,\infty},K_{3,\infty};\mathcal{W}(\itSigma_\infty,\psi_{3,\F_\infty}) \otimes M_{\mu,\C}^\vee\right)\\
& = \bigotimes_{v \in S_\infty}H^3(\frak{g}_3,K_3;\mathcal{W}(\itSigma_v,\psi_{3,\R}) \otimes M_{\mu_v,\C}^\vee).
\end{split}
\end{align}
For $\varphi \in \itSigma$, let $W_{\varphi,\psi_{3,\F}}$ be the Whittaker function of $\varphi$ with respect to $\psi_{3,\F}$ defined by
\[
W_{\varphi,\psi_{3,\F}}(g) = \int_{N_3(\F)\backslash N_3(\A_\F)}\varphi(ug)\overline{\psi_{3,\F}(u)}\,du^{\rm Tam}.
\]
Here $du^{\rm Tam}$ is the Haar measure on $N_3(\A_\F)$ such that ${\rm vol}(N_3(\F)\backslash N_3(\A_\F), du^{\rm Tam})=1$.
We have the $\GL_3(\A_{\F,f})$-equivariant isomorphism 
\begin{align*}
H^{2d}(\frak{g}_{3,\infty},K_{3,\infty};\itSigma \otimes M_{\mu,\C}^\vee) \longrightarrow \mathcal{W}(\itSigma_f,\psi_{3,\F}^{(\infty)}),\quad \varphi \longmapsto W_\varphi^{(\infty)}
\end{align*}
defined as follows: For $\varphi \in H^{2d}(\frak{g}_{3,\infty},K_{3,\infty};\itSigma \otimes M_{\mu,\C}^\vee)$, by Lemma \ref{L:cohomology generator} we can write
\begin{align}\label{E:2.3}
\varphi = \sum_{\underline{i}}\varphi_{\underline{i}}\otimes \bigotimes_{v \in S_\infty}\frac{(\sqrt{-1})^{{\sf w}/2}}{(\ell_v+i_v)!}\cdot E_-^{\ell_v+i_v}\cdot \left(X_{-1}^*\wedge X_{-2}^* \otimes \rho_{\mu_v^\vee}\left( \bp 1 & 0 & 1  \\ \sqrt{-1} & 0 & -\sqrt{-1} \\ 0 & 1 & 0\ep\right)P_{\mu_v^\vee}^+\right)
\end{align}
for some uniquely determined $\varphi_{\underline{i}} \in \itSigma$. 
Here $\underline{i} = \prod_{v \in S_\infty}i_v $ with $-\ell_v \leq i_v \leq \ell_v$.
Then $W_\varphi^{(\infty)} \in \mathcal{W}(\itSigma_f,\psi_{3,\F}^{(\infty)})$ is the unique Whittaker function depending only on $\varphi$ such that
\[
W_{\varphi_{\underline{i}},\psi_{3,\F}} = \prod_{v \in S_\infty}W_{(\ell_v, {\sf w}, \epsilon_v;\,i_v)}\cdot W_\varphi^{(\infty)}.
\]
For $\sigma \in {\rm Aut}(\C)$, we define the $\GL_3(\A_{\F,f})$-equivariant isomorphism 
\begin{align*}
H^{2d}(\frak{g}_{3,\infty},K_{3,\infty};{}^\sigma\!\itSigma \otimes M_{{}^\sigma\!\mu,\C}^\vee) \longrightarrow \mathcal{W}({}^\sigma\!\itSigma_f,\psi_{3,\F}^{(\infty)}),\quad \varphi \longmapsto W_\varphi^{(\infty)}
\end{align*}
in a similar way by considering the class 
\[
[{}^\sigma\!\itSigma_\infty]_b = \bigotimes_{v \in S_\infty}[\itSigma_{\sigma^{-1}\circ v}]_b.
\]
Let
\[
H^{2d}(\frak{g}_{3,\infty},K_{3,\infty};\itSigma \otimes M_{\mu,\C}^\vee) \longrightarrow  H^{2d}(\frak{g}_{3,\infty},K_{3,\infty};{}^\sigma\!\itSigma \otimes M_{\mu,\C}^\vee),\quad \varphi \longmapsto {}^\sigma\!\varphi
\]
be the $\sigma$-linear $\GL_3(\A_{\F,f})$-equivariant isomorphism such that the diagram 
\[
\begin{tikzcd}
H^{2d}(\frak{g}_{3,\infty},K_{3,\infty};\itSigma \otimes M_{\mu,\C}^\vee) \arrow[r] \arrow[d] & H^{2d}(\frak{g}_{3,\infty},K_{3,\infty};{}^\sigma\!\itSigma \otimes M_{{}^\sigma\!\mu,\C}^\vee) \arrow[d]\\
\mathcal{W}(\itSigma_f,\psi_{3,\F}^{(\infty)}) \arrow[r, "t_\sigma"]  & \mathcal{W}({}^\sigma\!\itSigma_f,\psi_{3,\F}^{(\infty)})
\end{tikzcd}
\]
is commute.
Similarly, replacing $(\sqrt{-1})^{{\sf w}/2}$ and $X_{-1}^* \wedge X_{-2}^*$ by $(\sqrt{-1})^{1+{\sf w}/2}$ and $X_0^*\wedge X_{-1}^* \wedge X_{-2}^*$ in (\ref{E:2.3}), we define the $\GL_3(\A_{\F,f})$-equivariant isomorphism
\begin{align*}\label{E:Galois equiv 2}
H^{3d}(\frak{g}_{3,\infty},K_{3,\infty};\itSigma \otimes M_{\mu,\C}^\vee) \longrightarrow \mathcal{W}(\itSigma_f,\psi_{3,\F}^{(\infty)}),\quad \varphi \longmapsto W_\varphi^{(\infty)}
\end{align*}
and the $\sigma$-linear $\GL_3(\A_{\F,f})$-equivariant isomorphism
\[
H^{3d}(\frak{g}_{3,\infty},K_{3,\infty};\itSigma \otimes M_{\mu,\C}^\vee) \longrightarrow  H^{3d}(\frak{g}_{3,\infty},K_{3,\infty};{}^\sigma\!\itSigma \otimes M_{{}^\sigma\!\mu,\C}^\vee), \quad \varphi \longmapsto {}^\sigma\!\varphi.
\]
Comparing the $\Q(\itSigma)$-rational structures given by the cuspidal cohomology and by the Whittaker model, we obtain the following lemma.
\begin{lemma}\label{L:Whittaker periods}
There exist $p^b(\itSigma), p^t(\itSigma) \in \C^\times$, unique up to $\Q(\itSigma)^\times$, such that
\begin{align*}
{\rm Image}\left(H^{2d}(\frak{g}_{3,\infty},K_{3,\infty};\itSigma \otimes M_{\mu,\C}^\vee)^{{\rm Aut}(\C/\Q(\itSigma))}\right)& = \frac{\mathcal{W}(\itSigma_f,\psi_{3,\F}^{(\infty)})^{{\rm Aut}(\C/\Q(\itSigma))}}{p^b(\itSigma)},\\
{\rm Image}\left(H^{3d}(\frak{g}_{3,\infty},K_{3,\infty};\itSigma \otimes M_{\mu,\C}^\vee)^{{\rm Aut}(\C/\Q(\itSigma))}\right) & = \frac{\mathcal{W}(\itSigma_f,\psi_{3,\F}^{(\infty)})^{{\rm Aut}(\C/\Q(\itSigma))}}{p^t(\itSigma)}.
\end{align*}
Moreover, we can normalize the periods so that
\[
T_\sigma\left(\frac{\varphi_1}{p^b(\itSigma)}\right) = \frac{{}^\sigma\!\varphi_1}{p^b({}^\sigma\!\itSigma)},\quad 
T_\sigma\left(\frac{\varphi_2}{p^t(\itSigma)}\right) = \frac{{}^\sigma\!\varphi_2}{p^t({}^\sigma\!\itSigma)}
\]
for all $\sigma\in {\rm Aut}(\C)$, $\varphi_1 \in H^{2d}(\frak{g}_{3,\infty},K_{3,\infty};\itSigma \otimes M_{\mu,\C}^\vee)$ and $\varphi_2 \in H^{3d}(\frak{g}_{3,\infty},K_{3,\infty};\itSigma \otimes M_{\mu,\C}^\vee)$.
\end{lemma}
We call $p^b(\itSigma)$ and $p^t(\itSigma)$ the bottom degree Whittaker period and the top degree Whittaker period, respectively, of $\itSigma$.

\section{Proof of Main Result}\label{S:proof}

Let $\itPi$ be a cohomological irreducible cuspidal automorphic representation of $\GL_2(\A_\F)$ with central character $\omega_\itPi$. We have $|\omega_\itPi| = |\mbox{ }|_{\A_\F}^{\sf w}$ for some ${\sf w} \in \Z$.
For $v \in S_\infty$, we have
\[
\itPi_v = D_{\kappa_v} \otimes |\mbox{ }|^{{\sf w}/2}
\]
for some $\kappa_v \in \Z_{\geq 2}$ such that $\kappa_v \equiv {\sf w} \,({\rm mod}\,2)$.
Recall that a newform of $\itPi$ is a cusp form in $\itPi$ which has weight $(\kappa_v)_{v \in S_\infty}$ under the $\prod_{v \in S_\infty}{\rm SO}(2)$-action and is new for the $\GL_2(\A_{\F,f})$-action in the sense of Casselman \cite{Casselman1973}. Newforms of $\itPi$ are unique up to scalars and we let $f_\itPi$ be the normalized newform so that
\begin{align*}
W_{f_\itPi,\psi_{2,\F}}({\rm diag}(d_\F^{-1},1)) = \int_{\F\backslash \A_\F}f_\itPi\left(\bp 1 & x \\ 0 & 1\ep {\rm diag}(d_\F^{-1},1)\right)\overline{\psi_{\F}(x)}\,dx^{\rm Tam} = e^{-2\pi d}.
\end{align*}
Here ${\rm vol}(\F\backslash \A_\F,dx^{\rm Tam})=1$ and $d_\F \in \A_{\F,f}^\times$ is any element such that $d_{\F,v} \in \F_v^\times$ generates the different ideal of $\F_v$ for each finite place $v$.
Similarly we define the normalized newform of $\itPi^\vee$. Let $\Vert f_\itPi \Vert$ be the Petersson norm of $f_\itPi$ defined by
\begin{align}\label{E:Petersson}
\Vert f_\itPi \Vert = \int_{\A_\F^\times\GL_2(\F)\backslash \GL_2(\A_\F)}f_\itPi(g)f_{\itPi^\vee}(g\cdot {\rm diag}(-1,1)_\infty)\,dg^{\rm Tam}.
\end{align}
Note that by definition we have $\Vert f_\itPi \Vert =\Vert f_{\itPi^\vee} \Vert$.
Let ${\rm Sym}^2\itPi$ be the functorial lifts of $\itPi$ with respect to the symmetric square representation of $\GL_2$. The functoriality was established by Gelbart--Jacquet \cite{GJ1978}.
The following theorem is our main result which expresses the algebraicity of the non-critical value $L^{(\infty)}(1,\itPi, {\rm Sym}^4 \otimes \omega_\itPi^{-2})$ in terms of the top degree Whittaker period $p^t({\rm Sym}^2\itPi)$ of ${\rm Sym}^2\itPi$ in Lemma \ref{L:Whittaker periods} and the Petersson norm $\Vert f_\itPi \Vert$.

\begin{thm}
Suppose $\itPi$ is non-CM and $\kappa_v \geq 3$ for all $v \in S_\infty$. We have
\begin{align*}
\sigma \left( \frac{L^{(\infty)}(1,\itPi,{\rm Sym}^4 \otimes \omega_\itPi^{-2})}{\pi^{3\sum_{v \in S_\infty}\kappa_v}\cdot G(\omega_\itPi)^{-3}\cdot \Vert f_\itPi\Vert\cdot p^t({\rm Sym}^2\itPi)}\right) = \frac{L^{(\infty)}(1,{}^\sigma\!\itPi, {\rm Sym}^4 \otimes {}^\sigma\!\omega_\itPi^{-2})}{\pi^{3\sum_{v \in S_\infty}\kappa_v}\cdot G({}^\sigma\!\omega_\itPi)^{-3}\cdot \Vert f_{{}^\sigma\!\itPi}\Vert\cdot p^t({\rm Sym}^2{}^\sigma\!\itPi)}
\end{align*}
for all $\sigma \in {\rm Aut}(\C)$. Here $G(\omega_\itPi)$ is the Gauss sum of $\omega_\itPi$ and $p^t({\rm Sym}^2\itPi)$ is the top degree Whittaker period of ${\rm Sym}^2\itPi$.
\end{thm}

\begin{proof}
We have the factorization of $L$-functions
\[
L(s,{\rm Sym}^2\itPi \times {\rm Sym}^2\itPi) = L(s,\itPi , {\rm Sym}^4\otimes \omega_\itPi^{-2})\cdot L(s,\itPi , {\rm Sym}^2\otimes \omega_\itPi^{-1}) \cdot  \zeta_\F(s).
\]
Here $\zeta_\F(s)$ is the completed Dedekind zeta function of $\F$.
By the result of Im \cite[Corollary 5.4]{Im1991} (see also Sturm \cite{Sturm1989}), we have
\begin{align*}
\sigma \left( \frac{L^{(\infty)}(1,\itPi,{\rm Sym}^2 \otimes \omega_\itPi^{-1})}{\pi^{2d+\sum_{v \in S_\infty}\kappa_v}\cdot \Vert f_\itPi \Vert} \right) = \frac{L^{(\infty)}(1,{}^\sigma\!\itPi,{\rm Sym}^2 \otimes {}^\sigma\!\omega_\itPi^{-1})}{\pi^{2d+\sum_{v \in S_\infty}\kappa_v}\cdot \Vert f_{{}^\sigma\!\itPi} \Vert}.
\end{align*}
for all $\sigma \in {\rm Aut}(\C)$.
By the class number formula, we have
\[
{\rm Res}_{s=1}\zeta_\F^{(\infty)}(s) \in |D_\F|^{1/2}\cdot {\rm Reg}_\F\cdot \Q^\times.
\]
Note that ${\rm Sym}^2 \itPi$ contributes to the cuspidal cohomology of $\GL_3(\A_\F)$ with coefficients in $M^\vee_{\mu,\C}$ with $\mu = \prod_{v \in S_\infty}\mu_v$ and 
\[
\mu_v = (\kappa_v-2+{\sf w},{\sf w},-\kappa_v+2+{\sf w})
\]
for $v \in S_\infty$.
By Theorems \ref{T:GJ period relation} and \ref{T:adjoint} below, we have
\begin{align*}
\sigma\left(\frac{|D_\F|^{1/2}\cdot \pi^{2d}\cdot G(\omega_\itPi)^{-3}\cdot\Vert f_\itPi \Vert^2}{p^b({\rm Sym}^2\itPi^\vee)}\right) &= \frac{|D_\F|^{1/2}\cdot \pi^{2d}\cdot G({}^\sigma\!\omega_\itPi)^{-3}\cdot\Vert f_{{}^\sigma\!\itPi} \Vert^2}{p^b({\rm Sym}^2 {}^\sigma\!\itPi^\vee)},\\
\sigma \left(\frac{{\rm Res}_{s=1}L(s,{\rm Sym}^2\itPi \times {\rm Sym}^2\itPi)}{{\rm Reg}_\F\cdot\pi^{4\sum_{v \in S_\infty}\kappa_v}\cdot p^t({\rm Sym}^2\itPi)\cdot p^b({\rm Sym}^2\itPi^\vee)}\right) &= \frac{{\rm Res}_{s=1}L(s,{\rm Sym}^2{}^\sigma\!\itPi \times {\rm Sym}^2{}^\sigma\!\itPi)}{{\rm Reg}_\F\cdot\pi^{4\sum_{v \in S_\infty}\kappa_v}\cdot p^t({\rm Sym}^2{}^\sigma\!\itPi)\cdot p^b({\rm Sym}^2{}^\sigma\!\itPi^\vee)} 
\end{align*}
respectively, for all $\sigma \in {\rm Aut}(\C)$.
Here $p^b({\rm Sym}^2\itPi^\vee)$ is the bottom degree Whittaker period of ${\rm Sym}^2\itPi^\vee$.
The assertion thus follows. This completes the proof.
\end{proof}

\section{Algebraicity of the Rankin--Selberg $L$-functions for $\GL_3 \times \GL_2$}\label{S:GL_3 GL_2}

Let $\itPi = \bigotimes_v \itPi_v$ and $\itSigma = \bigotimes_v \itSigma_v$ be cohomological irreducible cuspidal automorphic representations of $\GL_2(\A_\F)$ and $\GL_3(\A_\F)$ with central characters $\omega_\itPi$ and $\omega_\itSigma$, respectively.
We have $|\omega_\itPi|=|\mbox{ }|_{\A_\F}^{{\sf w}(\itPi)}$ and $|\omega_\itSigma|=|\mbox{ }|_{\A_\F}^{{\sf w}(\itSigma)}$ for some ${\sf w}(\itPi) \in \Z$ and even integer ${\sf w}(\itSigma)$.
Let 
\[
\lambda = \prod_{v \in S_\infty}\lambda_v \in \prod_{v \in S_\infty}X^+(T_2),\quad \mu =\prod_{v \in S_\infty} \mu_v \in \prod_{v \in S_\infty}X^+(T_3)
\] such that $\itPi$ and $\itSigma$ contribute to the cuspidal cohomology of $\GL_2(\A_\F)$ and $\GL_3(\A_\F)$ with coefficients in $M_{\lambda,\C}^\vee$ and $M_{\mu,\C}^\vee$, respectively.
For $v \in S_\infty$, we have
\[
\lambda_v = (\tfrac{\kappa_v-2+{\sf w}(\itPi)}{2}, \tfrac{-\kappa_v+2+{\sf w}(\itPi)}{2}),\quad \mu_v = (\tfrac{\ell_v-3+{\sf w}(\itSigma)}{2},\tfrac{{\sf w}(\itSigma)}{2},\tfrac{-\ell_v+3+{\sf w}(\itSigma)}{2})
\]
for some $\kappa_v \in \Z_{\geq 2}$ with $\kappa_v \equiv {\sf w}(\itPi)\,({\rm mod }\,2)$, odd integer $\ell_v \in \Z_{\geq 3}$, and
\[
\itPi_v = D_{\kappa_v} \otimes |\mbox{ }|^{{\sf w}(\itPi)/2},\quad \itSigma_v = {\rm Ind}^{\GL_3(\R)}_{P_{2,1}(\R)}(D_{\ell_v} \boxtimes \epsilon_v) \otimes |\mbox{ }|^{{\sf w}(\itSigma)/2}
\]
for some quadratic character $\epsilon_v$ of $\R^\times$.
Let 
\[
L(s,\itSigma \times \itPi)
\]
be the Rankin--Selberg $L$-function of $\itSigma \times \itPi$. We denote by $L^{(\infty)}(s,\itSigma \times \itPi)$ the $L$-function obtained by excluding the archimedean $L$-factors. A critical point for $L(s,\itSigma \times \itPi)$ is a half-integer $m+\tfrac{1}{2}$ which is not a pole of the archimedean local factors $L(s,\itSigma_v \times \itPi_v)$ and $L(1-s,\itSigma_v^\vee \times \itPi_v^\vee)$ for all $v \in S_\infty$.
The aim of this section is to prove Theorem \ref{T:RS} on the algebraicity of the critical values of $L(s,\itSigma \times \itPi)$ in terms of the bottom degree Whittaker periods of $\itSigma$ and $\itPi$ under the assumption that $\ell_v > \kappa_v$ for all $v \in S_\infty$. As a consequence, we obtain in Theorem \ref{T:GJ period relation} a period relation between the bottom degree Whittaker period of the Gelbart--Jacquet lift of $\itPi$ and the Petersson norm of the normalized newform of $\itPi$. In \S\,\ref{SS:Deligne}, we propose in Conjecture \ref{C:Deligne-Whittaker} a conjectural period relation between the bottom degree Whittaker period of $\itSigma$ and the product of Deligne's motivic periods for the conjectural motive attached to $\itSigma$.
We will show that the conjecture is a consequence of the compatibility of Theorem \ref{T:RS} with Deligne's conjecture \cite{Deligne1979} for $\GL_3 \times \GL_2$.

\subsection{Algebraicity for $\GL_2 \times \GL_1$}

In this section, we recall the result of Raghuram--Tanabe 
\cite{RT2011} on the algebraicity of the critical values of the twisted standard $L$-functions of $\itPi$ in terms of the Whittaker periods of $\itPi$.
Firstly we recall the definition of the Whittaker periods of $\itPi$.
Recall 
\[
(\rho_{\lambda^\vee}, M_{\lambda}^\vee) = (\otimes_{v \in S_\infty}\rho_{\lambda_v^\vee},\otimes_{v \in S_\infty}M_{\lambda_v,\C}^\vee).
\] 
Similarly as in \S\,\ref{SS:cohomology}, we define the cuspidal cohomology
\[
H_{\rm cusp}^d(\mathcal{S}_2,\mathcal{M}^\vee_{\lambda,\C}) = H^1(\frak{g}_{2,\infty},K_{2,\infty};\mathcal{A}_{\rm cusp}(\GL_2(\F)\backslash \GL_2(\A_\F)) \otimes M_{\lambda,\C}^\vee)
\]
of $\GL_2(\A_\F)$ with coefficients in $M_{\lambda,\C}^\vee$.
Then $\itPi_f = \bigotimes_{v \nmid \infty}\itPi_v$ appears in $H_{\rm cusp}^d(\mathcal{S}_2,\mathcal{M}^\vee_{\lambda,\C})$ with multiplicity $2^d$ and the $\itPi_f$-isotypic component is equal to 
\[
H^d(\frak{g}_{2,\infty},K_{2,\infty};\itPi \otimes M_{\lambda,\C}^\vee) = \left(\itPi \otimes  (\frak{g}_{2,\infty,\C} / \frak{k}_{2,\infty,\C})^* \otimes M_{\lambda,\C}^\vee\right)^{K_{2,\infty}}.
\]
Note that $\prod_{v \in S_\infty}{\rm O}(2) / {\rm SO}(2)$ acts on $H^d(\frak{g}_{2,\infty},K_{2,\infty};\itPi \otimes M_{\lambda,\C}^\vee)$ in a natural way.
For $\underline{\varepsilon} = (\varepsilon_v)_{v \in S_\infty} \in \{\pm1\}^{S_\infty}$, let
\[
H^d(\frak{g}_{2,\infty},K_{2,\infty};\itPi \otimes M_{\lambda,\C}^\vee)[\underline{\varepsilon}] 
\]
be the $1$-dimensional subspace of $H^d(\frak{g}_{2,\infty},K_{2,\infty};\itPi \otimes M_{\lambda,\C}^\vee)$ such that ${\rm O}(2) / {\rm SO}(2)$ acts by $\varepsilon_v$ for each $v \in S_\infty$.
Let $v \in S_\infty$.
Similarly we denote by
\[
H^1(\frak{g}_2,K_2;\mathcal{W}(\itPi_v,\psi_{2,\R}) \otimes M_{\lambda_v,\C}^\vee)[\pm] = \left(\mathcal{W}(\itPi_v,\psi_{2,\R}) \otimes  (\frak{g}_{2,\C} / \frak{k}_{2,\C})^* \otimes M_{\lambda_v,\C}^\vee\right)^{K_2,\pm}
\]
the $\pm 1$-eigenspace under the ${\rm O}(2) / {\rm SO}(2)$-action.
Let $W_{\kappa_v}^\pm$ be the Whittaker function of $D_{\kappa_v}$ with respect to $\psi_{2,\R}$ of weight $\pm \kappa_v$ normalized so that
\begin{align}\label{E:GL_2 Whittaker}
W_{\kappa_v}^\pm ({\rm diag}(a,1)) = |a|^{\kappa_v/2}e^{-2\pi|a|}\cdot\mathbb{I}_{\R_+}(\pm a).
\end{align}
Following the normalization in \cite[(3.6) and (3.7)]{RT2011}, we define the class 
\[
[\itPi_v]^\pm\in  H^1(\frak{g}_2,K_2;\mathcal{W}(\itPi_v,\psi_{2,\R}) \otimes M_{\lambda_v,\C}^\vee)[\pm]
\]
by
\begin{align*}
[\itPi_v]^\pm &= (W_{\kappa_v}^+\otimes|\mbox{ }|^{{\sf w}(\itPi)/2}) \otimes Y_+^*\otimes (\sqrt{-1}\,x+y)^{\kappa_v-2}\\
&\pm (\sqrt{-1})^{{\sf w}(\itPi)}(W_{\kappa_v}^-\otimes|\mbox{ }|^{{\sf w}(\itPi)/2}) \otimes Y_-^*\otimes (x+\sqrt{-1}\,y)^{\kappa_v-2}.
\end{align*} 
For $\underline{\varepsilon} = (\varepsilon_v)_{v \in S_\infty} \in \{\pm1\}^{S_\infty}$, let
\begin{align}\label{E:GL_2 class}
\begin{split}
[\itPi_\infty]^{\underline{\varepsilon}} = \bigotimes_{v \in S_\infty}[\itPi_v]^{\varepsilon_v} &\in H^d(\frak{g}_{2,\infty},K_{2,\infty};\mathcal{W}(\itPi_\infty,\psi_{2,\F_\infty})\otimes M_{\lambda,\C}^\vee)[\underline{\varepsilon}]\\
& =\bigotimes_{v \in S_\infty} H^1(\frak{g}_2,K_2;\mathcal{W}(\itPi_v,\psi_{2,\R})\otimes M_{\lambda_v,\C}^\vee)[\varepsilon_v]
\end{split}
\end{align}
We have the $\GL_2(\A_{\F,f})$-equivariant isomorphism 
\begin{align*}
H^d(\frak{g}_{2,\infty},K_{2,\infty};\itPi \otimes M_{\lambda,\C}^\vee)[\underline{\varepsilon}] \longrightarrow \mathcal{W}(\itPi_f,\psi_{2,\F}^{(\infty)}),\quad \varphi \longmapsto W_\varphi^{(\infty)}
\end{align*}
defined as follows:
For $\varphi \in H^d(\frak{g}_{2,\infty},K_{2,\infty};\itPi \otimes M_{\lambda,\C}^\vee)[\underline{\varepsilon}]$, we can write
\begin{align*}
\varphi = \sum_{I \subseteq S_\infty}\varphi_{I}\otimes\bigotimes_{v \in S_\infty,\,v\in I}Y_+^* \otimes (\sqrt{-1}\,x+y)^{\kappa_v-2} \otimes \bigotimes_{v \in S_\infty,\,v\notin I}\varepsilon_v\cdot (\sqrt{-1})^{{\sf w}(\itPi)}\cdot Y_-^*\otimes (x+\sqrt{-1}\,y)^{\kappa_v-2}.
\end{align*}
Let $W_\varphi^{(\infty)} \in \mathcal{W}(\itPi_f,\psi_{2,\F}^{(\infty)})$ be the unique Whittaker function such that 
\[
W_{\varphi_{I},\psi_{2,\F}} = \prod_{v \in S_\infty,\,v \in I}(W_{\kappa_v}^+\otimes|\mbox{ }|^{{\sf w}(\itPi)/2}) \prod_{v \in S_\infty,\,v \notin I}(W_{\kappa_v}^-\otimes|\mbox{ }|^{{\sf w}(\itPi)/2}) \cdot W_\varphi^{(\infty)}
\]
for all $I \subseteq S_\infty$.
For $\sigma \in {\rm Aut}(\C)$, similarly we define the $\GL_2(\A_{\F,f})$-equivariant isomorphism 
\begin{align*}
H^d(\frak{g}_{2,\infty},K_{2,\infty};{}^\sigma\!\itPi \otimes M_{{}^\sigma\!\lambda,\C}^\vee)[{}^\sigma\!\underline{\varepsilon}] \longrightarrow \mathcal{W}(\itPi_f,\psi_{2,\F}^{(\infty)}),\quad \varphi \longmapsto W_\varphi^{(\infty)}
\end{align*}
by considering the class
\[
[{}^\sigma\!\itPi_\infty]^{{}^\sigma\!\underline{\varepsilon}} = \bigotimes_{v \in S_\infty} [\itPi_{\sigma^{-1}\circ v}]^{\varepsilon_{\sigma^{-1}\circ v}}.
\]
Here ${}^\sigma\!\underline{\varepsilon} = (\varepsilon_{\sigma^{-1}\circ v})_{v \in S_\infty}$.
Similarly as in (\ref{E:sigma linear map 1}) and (\ref{E:sigma linear map 2}), we have the $\sigma$-linear $\GL_2(\A_{\F,f})$-equivariant isomorphisms
\begin{align*}
T_\sigma : H^d(\frak{g}_{2,\infty},K_{2,\infty};\itPi \otimes M_{\lambda,\C}^\vee)[\underline{\varepsilon}] &\longrightarrow H^d(\frak{g}_{2,\infty},K_{2,\infty};\itPi \otimes M_{{}^\sigma\!\lambda,\C}^\vee)[{}^\sigma\!\underline{\varepsilon}],\\
t_\sigma : \mathcal{W}(\itPi_f,\psi_{2,\F}^{(\infty)}) & \longrightarrow \mathcal{W}({}^\sigma\!\itPi_f,\psi_{2,\F}^{(\infty)}).
\end{align*}
The Whittaker period $p(\itPi,\underline{\varepsilon}) \in \C^\times$ is the non-zero complex number, unique up to $\Q(\itPi)^\times$, such that
\[
{\rm Image}\left(H^d(\frak{g}_{2,\infty},K_{2,\infty};\itPi \otimes M_{\lambda,\C}^\vee)[\underline{\varepsilon}]^{{\rm Aut}(\C/\Q(\itPi))}\right) = \frac{\mathcal{W}(\itPi_f,\psi_{2,\F}^{(\infty)})^{{\rm Aut}(\C/\Q(\itPi))}}{p(\itPi,\underline{\varepsilon})}.
\]
Let 
\[
H^d(\frak{g}_{2,\infty},K_{2,\infty};\itPi \otimes M_{\lambda,\C}^\vee)[\underline{\varepsilon}] \longrightarrow  H^d(\frak{g}_{2,\infty},K_{2,\infty};{}^\sigma\!\itPi \otimes M_{{}^\sigma\!\lambda,\C}^\vee)[{}^\sigma\!\underline{\varepsilon}],\quad \varphi \longmapsto {}^\sigma\!\varphi
\]
be the $\sigma$-linear $\GL_2(\A_{\F,f})$-equivariant isomorphism such that the diagram 
\[
\begin{tikzcd}
H^d(\frak{g}_{2,\infty},K_{2,\infty};\itPi \otimes M_{\lambda,\C}^\vee)[\underline{\varepsilon}] \arrow[r] \arrow[d] & H^d(\frak{g}_{2,\infty},K_{2,\infty};{}^\sigma\!\itPi \otimes M_{{}^\sigma\!\lambda,\C}^\vee)[{}^\sigma\!\underline{\varepsilon}] \arrow[d]\\
\mathcal{W}(\itPi_f,\psi_{2,\F}^{(\infty)}) \arrow[r, "t_\sigma"]  & \mathcal{W}({}^\sigma\!\itPi_f,\psi_{2,\F}^{(\infty)})
\end{tikzcd}
\]
is commute.
Then we can normalize the Whittaker periods so that
\[
T_\sigma\left(\frac{\varphi}{p(\itPi,\underline{\varepsilon})}\right) = \frac{{}^\sigma\!\varphi}{p({}^\sigma\!\itPi,{}^\sigma\!\underline{\varepsilon})}
\]
for all $\sigma \in {\rm Aut}(\C)$ and $\varphi \in H^d(\frak{g}_{2,\infty},K_{2,\infty};\itPi \otimes M_{\lambda,\C}^\vee)[\underline{\varepsilon}]$.

We prove in the following lemma a period relation between the product of the Whittaker periods of $\itPi$ and the Petersson norm of the normalized newform of $\itPi$. 
It should be well-known and is a consequence of the Poincar\'e duality. We include a proof here for the reader's convenience.

\begin{lemma}\label{L:period relation GL_2}
Let $\underline{\varepsilon}\in \{\pm1\}^{S_\infty}$.
We have
\[
\sigma\left(\frac{p(\itPi,\underline{\varepsilon})\cdot p(\itPi,-\underline{\varepsilon})}{(2\pi\sqrt{-1})^d\cdot(\sqrt{-1})^{d{\sf w}(\itPi)}\cdot G(\omega_\itPi)\cdot\Vert f_\itPi \Vert}\right) = \frac{p({}^\sigma\!\itPi,{}^\sigma\!\underline{\varepsilon})\cdot p({}^\sigma\!\itPi,-{}^\sigma\!\underline{\varepsilon})}{(2\pi\sqrt{-1})^d\cdot(\sqrt{-1})^{d{\sf w}(\itPi)}\cdot G({}^\sigma\!\omega_\itPi)\cdot\Vert f_{{}^\sigma\!\itPi} \Vert}
\]
for all $\sigma \in {\rm Aut}(\C)$. Here $\Vert f_\itPi \Vert$ is the Petersson norm of $f_\itPi$ in (\ref{E:Petersson}).
\end{lemma}

\begin{proof}

Let $H^d_c(\mathcal{S}_2,\mathcal{M}_{\lambda,\C}^\vee)$ be the cohomology with compact support. We have the Poincar\'e pairing
\[
H^d_c(\mathcal{S}_2,\mathcal{M}_{\lambda,\C}^\vee) \otimes H^d(\mathcal{S}_2,\mathcal{M}_{\lambda,\C}) \longrightarrow \C,\quad c_1 \otimes c_2 \longmapsto \int_{\mathcal{S}_2}c_1\wedge c_2
\]
which satisfies the Galois-equivariant property
\[
\sigma \left( \int_{\mathcal{S}_2}c_1\wedge c_2\right) = \int_{\mathcal{S}_2}T_\sigma c_1\wedge T_\sigma c_2
\]
for all $\sigma \in {\rm Aut}(\C)$ (cf.\,\cite[p.\,124]{Clozel1990}).
For $I \subseteq S_\infty$, define $f_{\itPi,I} \in \itPi$ by
\[
f_{\itPi,I}(g) = f_{\itPi}\left(g\cdot \prod_{v \in S_\infty,v \notin I}{\rm diag}(-1,1)\right).
\]
Let 
\[
\varphi_\itPi \in H^d(\frak{g}_{2,\infty},K_{2,\infty};\itPi \otimes M_{\lambda,\C}^\vee)[\underline{\varepsilon}],\quad \varphi_{\itPi^\vee} \in H^d(\frak{g}_{2,\infty},K_{2,\infty};\itPi^\vee \otimes M_{\lambda,\C})[-\underline{\varepsilon}]
\]
be the classes defined by
\begin{align*}
\varphi_\itPi &= \sum_{I \subseteq S_\infty}f_{\itPi,I}\otimes\bigotimes_{v \in S_\infty,\,v\in I}Y_+^* \otimes (\sqrt{-1}\,x+y)^{\kappa_v-2} \otimes \bigotimes_{v \in S_\infty,\,v\notin I}\varepsilon_v\cdot (\sqrt{-1})^{{\sf w}(\itPi)}\cdot Y_-^*\otimes (x+\sqrt{-1}\,y)^{\kappa_v-2},\\
\varphi_{\itPi^\vee} &= \sum_{J \subseteq S_\infty}f_{\itPi^\vee,J}\otimes\bigotimes_{v \in S_\infty,\,v\in J}Y_+^* \otimes (\sqrt{-1}\,x+y)^{\kappa_v-2} \otimes \bigotimes_{v \in S_\infty,\,v\notin J}(-\varepsilon_v)\cdot (\sqrt{-1})^{-{\sf w}(\itPi)}\cdot Y_-^*\otimes (x+\sqrt{-1}\,y)^{\kappa_v-2}.
\end{align*}
By the definition of $f_\itPi$ and $f_{\itPi^\vee}$, it is easy to see that
\[
{}^\sigma \varphi_\itPi = \varphi_{{}^\sigma\!\itPi},\quad {}^\sigma \varphi_{\itPi^\vee} = \varphi_{{}^\sigma\!\itPi^\vee}
\]
for all $\sigma \in {\rm Aut}(\C)$.
Note that
\begin{align*}
\int_{\A_\F^\times\GL_2(\F)\backslash\GL_2(\A_\F)}f_{
\itPi,I}(g)f_{\itPi^\vee,J}(g)\,dg^{\rm Tam} = \begin{cases}
\Vert f_\itPi \Vert & \mbox{ if $J=S_\infty \setminus I$},\\
0 & \mbox{ otherwise}.
\end{cases}
\end{align*}
Also note that the comparison between the standatd measure $dg = \prod_v dg_v$ in \S\,\ref{SS:measure} and the Tamagawa measure on $\A_\F^\times \backslash \GL_2(\A_\F)$ is given by (cf.\,\cite[Lemma 6.1]{IP2018})
\[
dg = 2^{-d}\cdot|D_\F|^{3/2}\cdot \zeta_\F(2)\cdot dg^{\rm Tam},
\]
where $\zeta_\F(s)$ is the completed Dedekind zeta function of $\F$.
Now we evaluate the Poincar\'e pairing at $\varphi_\itPi \otimes \varphi_{\itPi^\vee}$ and obtain
\begin{align*}
&\int_{\mathcal{S}_2}\varphi_\itPi \wedge \varphi_{\itPi^\vee}\\
 &= C \cdot |D_\F|^{3/2}\cdot \zeta_\F(2)\cdot (\sqrt{-1})^{d}\cdot \Vert f_\itPi \Vert\\
&\times \prod_{v \in S_\infty}\varepsilon_v\cdot(\sqrt{-1})^{{\sf w}(\itPi)}\\
&\times\left(\left\<(x+\sqrt{-1}\,y)^{\kappa_v-2},(\sqrt{-1}\,x+y)^{\kappa_v-2}\right\>_{\lambda_v,\C}+(-1)^{{\sf w}(\itPi)} \left\<(\sqrt{-1}\,x+y)^{\kappa_v-2}, (x+\sqrt{-1}\,y)^{\kappa_v-2}\right\>_{\lambda_v,\C}\right)
\end{align*}
for some non-zero rational number $C$ independent of $\itPi$ and some non-zero $\GL_2(\C)$-equivariant bilinear pairing $\<\cdot,\cdot\>_{\lambda_v,\C}$ on $M_{\lambda_v,\C}^\vee \times M_{\lambda_v,\C}$ defined over $\Q$ for each $v \in S_\infty$.
Here the extra factor $(\sqrt{-1})^d$ in the first line is due to 
\[
Y_+^* \wedge Y_-^* =  \bp1 & 0 \\ 0 & 0\ep^* \wedge \bp 0 & 1 \\ 0 & 0\ep^* \otimes 8 \sqrt{-1}
\]
in $\bigwedge^2 (\frak{g}_{2,\C} / \frak{k}_{2,\C})^*$  and $\bp1 & 0 \\ 0 & 0\ep^* \wedge \bp 0 & 1 \\ 0 & 0\ep^*$ correspondes to the standard measure on $\GL_2(\R) / K_2$.
It is easy to see that
\begin{align*}
&\left\<(x+\sqrt{-1}\,y)^{\kappa_v-2},(\sqrt{-1}\,x+y)^{\kappa_v-2}\right\>_{\lambda_v,\C}+(-1)^{{\sf w}(\itPi)} \left\<(\sqrt{-1}\,x+y)^{\kappa_v-2}, (x+\sqrt{-1}\,y)^{\kappa_v-2}\right\>_{\lambda_v,\C}\\
& = 2^{\kappa_v-1}\cdot \<x^{\kappa_v-2},y^{\kappa_v-2}\>_{\lambda_v,\C} \in \Q^\times.
\end{align*}
We conclude that
\begin{align*}
&\sigma\left( \frac{|D_\F|^{3/2}\cdot \zeta_\F(2)\cdot (\sqrt{-1})^{d+d{\sf w}(\itPi)}\cdot \Vert f_\itPi \Vert}{p(\itPi,\underline{\varepsilon})\cdot p(\itPi^\vee,-\underline{\varepsilon})}\right)\\
 &=  C^{-1}\cdot 2^{d-\sum_{v \in S_\infty}\kappa_v}\prod_{v \in S_\infty}\varepsilon_v\cdot\<x^{\kappa_v-2},y^{\kappa_v-2}\>_{\lambda_v,\C}^{-1}\cdot \sigma \left( \int_{\mathcal{S}_2}\frac{\varphi_\itPi}{p(\itPi,\underline{\varepsilon})}\wedge \frac{\varphi_{\itPi^\vee}}{p(\itPi^\vee,-\underline{\varepsilon})}\right)\\
& = C^{-1}\cdot 2^{d-\sum_{v \in S_\infty}\kappa_v}\prod_{v \in S_\infty}\varepsilon_v\cdot\<x^{\kappa_v-2},y^{\kappa_v-2}\>_{\lambda_v,\C}^{-1}\cdot\int_{\mathcal{S}_2}\frac{\varphi_{{}^\sigma\!\itPi}}{p({}^\sigma\!\itPi,{}^\sigma\!\underline{\varepsilon})}\wedge \frac{\varphi_{{}^\sigma\!\itPi^\vee}}{p({}^\sigma\!\itPi^\vee,-{}^\sigma\!\underline{\varepsilon})}\\
& = \frac{|D_\F|^{3/2}\cdot \zeta_\F(2)\cdot (\sqrt{-1})^{d+d{\sf w}(\itPi)}\cdot \Vert f_{{}^\sigma\!\itPi} \Vert}{p({}^\sigma\!\itPi,{}^\sigma\!\underline{\varepsilon})\cdot p({}^\sigma\!\itPi^\vee,-{}^\sigma\!\underline{\varepsilon})}
\end{align*}
for all $\sigma \in {\rm Aut}(\C)$.
The assertion then follows from the period relation \cite{RS2008}
\[
\sigma \left( \frac{G(\omega_\itPi)\cdot p(\itPi^\vee,\pm\underline{\varepsilon})}{p(\itPi,\pm\underline{\varepsilon})}\right) = \frac{G({}^\sigma\!\omega_\itPi)\cdot p({}^\sigma\!\itPi^\vee,\pm{}^\sigma\!\underline{\varepsilon})}{p({}^\sigma\!\itPi,\pm{}^\sigma\!\underline{\varepsilon})}
\]
for all $\sigma \in {\rm Aut}(\C)$ and the fact that $|D_\F|^{3/2}\cdot \zeta_\F(2) \in \pi^d\cdot\Q^\times$ (cf.\,\cite[Proposition 3.1]{Shimura1978}).
This completes the proof.
\end{proof}

\begin{rmk}
The lemma can also be proved by combining the result \cite[Theorem 4.3-(II)]{Shimura1978} of Shimura and the period relation established by Raghuram--Tanabe \cite[(4.23)]{RT2011} between the Whittaker periods and the periods defined by Shimura \cite{Shimura1978}.
\end{rmk}

Let $\chi$ be a finite order Hecke character of $\A_\F^\times$ and 
\[
L(s,\itPi \times \chi)
\]
be the twisted standard $L$-function of $\itPi$ by $\chi$.
We denote by $L^{(\infty)}(s,\itPi \times \chi)$ the $L$-function obtained by excluding the archimedean $L$-factors.
The critical points for $L(s,\itPi\times\chi)$ are given by $m+\tfrac{1}{2} \in \Z+\tfrac{1}{2}$ with 
\[
\tfrac{-\min_{v \in S_\infty}\{\kappa_v\}-{\sf w}(\itPi)+2}{2} \leq m \leq \tfrac{\min_{v \in S_\infty}\{\kappa_v\}-{\sf w}(\itPi)-2}{2}.
\]
We have the following theorem of Raghuram--Tanabe \cite{RT2011}.

\begin{thm}[Raghuram--Tanabe]\label{T:RT2011}
Let $m+\tfrac{1}{2}$ be critical for $L(s,\itPi \times \chi)$. We have
\begin{align*}
&\sigma\left( \frac{L^{(\infty)}(m+\tfrac{1}{2},\itPi \times \chi)}{(2\pi\sqrt{-1})^{dm+\sum_{v \in S_\infty}(\kappa_v+{\sf w}(\itPi))/2}\cdot G(\chi) \cdot p(\itPi,(-1)^m\cdot{\rm sgn}(\chi))}\right)\\
&= \frac{L^{(\infty)}(m+\tfrac{1}{2},{}^\sigma\!\itPi \times {}^\sigma\!\chi)}{(2\pi\sqrt{-1})^{dm+\sum_{v \in S_\infty}(\kappa_v+{\sf w}(\itPi))/2}\cdot G({}^\sigma\!\chi) \cdot p({}^\sigma\!\itPi,(-1)^m\cdot{\rm sgn}(\chi))}
\end{align*}
for all $\sigma\in{\rm Aut}(\C)$.
\end{thm}

\begin{rmk}
Similar result was proved by Shimura \cite[Theorem 4.3-(I)]{Shimura1978} with the Whittaker periods replaced by the periods defined therein.
\end{rmk}

\subsection{Algebraicity for $\GL_3 \times \GL_2$}
Let $\iota : \GL_2 \rightarrow \GL_3$ be the embedding defined by $\iota(g) = \bp g & 0 \\ 0 & 1\ep$. 
Let $s: \bigwedge^2 (\frak{g}_{3,\C} / \frak{k}_{3,\C})^* \times (\frak{g}_{2,\C} / \frak{k}_{2,\C})^* \rightarrow \C$ be the bilinear homomorphism defined by
\[
\iota^*X^* \wedge \iota^* Y^* \wedge {\rm pr}(Z^*) = s(X^* \wedge Y^*, Z^*)\cdot \bp 1 & 0 \\ 0 & 1\ep^*\wedge \bp 1 & 0 \\ 0 & 0\ep^* \wedge \bp 0 & 1 \\ 0 & 0\ep^*.
\]
Here ${\rm pr}: (\frak{g}_{2,\C} / \frak{k}_{2,\C})^* \rightarrow (\frak{g}_{2,\C} / \frak{so}(2)_\C)^*$ is the natural surjection and $\iota^* : (\frak{g}_{3,\C} / \frak{k}_{3,\C})^*  \rightarrow (\frak{g}_{2,\C} / \frak{so}(2)_\C)^*$ is the homomorphism induced by the embedding $\iota$.
Let $v \in S_\infty$. For $W_1 \in \mathcal{W}(\itSigma_v, \psi_{3,\R})$ and $W_2 \in \mathcal{W}(\itPi_v, \psi_{2,\R})$, we define the local zeta integral
\[
Z_v(s,W_1,W_2) = \int_{N_2(\R)\backslash \GL_2(\R)}W_1(\iota(g))W_2({\rm diag}(-1,1)g)|\det(g)|^{s-1/2}\,dg.
\]
The integral converges absolutely for ${\rm Re}(s)$ sufficiently large and admits meromorphic continuation to $s \in \C$.
Moreover, the ratio
\[
\frac{Z_v(s,W_1,W_2)}{L(s,\itSigma_v\times \itPi_v)}
\]
is entire (cf.\,\cite{JS1990} and \cite{Jacquet2009}).
Let $m+\tfrac{1}{2} \in \Z+\tfrac{1}{2}$. 
Under the assumption $\ell_v > \kappa_v$, we know that $m+\tfrac{1}{2}$ is is not a pole of neither $L(s,\itSigma_v \times \itPi_v)$ nor $L(1-s,\itSigma_v^\vee \times \itPi_v^\vee)$ if and only if 
\[
{\rm Hom}_{\GL_2}(M_{\lambda_v+m}^\vee,M_{\mu_v})\neq 0.
\]
In this case, we fix a non-zero $\iota_{m,v} \in {\rm Hom}_{\GL_2}(M_{\lambda_v+m}^\vee,M_{\mu_v})$. We also fix a non-zero $\GL_3(\C)$-equivariant bilinear pairing $\<\cdot,\cdot\>_{\mu_v,\C} : M_{\mu_v,\C}^\vee \times M_{\mu_v,\C} \rightarrow \C$ defined over $\Q$.
Let 
\[
\<\cdot,\cdot\>_{\mu_v,\lambda_v,m}: \left(\bigwedge^2 (\frak{g}_{3,\C} / \frak{k}_{3,\C})^* \otimes M_{\mu_v,\C}^\vee\right)\times\left((\frak{g}_{2,\C} / \frak{k}_{2,\C})^* \otimes M_{\lambda_v,\C}^\vee\right) \longrightarrow \C
\]
and
\[
\<\cdot,\cdot\>_{m,v}: \left(\mathcal{W}(\itSigma_v,\psi_{3,\R}) \otimes \bigwedge^2 (\frak{g}_{3,\C} / \frak{k}_{3,\C})^* \otimes M_{\mu_v,\C}^\vee\right)\times\left(\mathcal{W}(\itPi_v,\psi_{2,\R}) \otimes  (\frak{g}_{2,\C} / \frak{k}_{2,\C})^* \otimes M_{\lambda_v,\C}^\vee\right) \longrightarrow \C
\]
be the bilinear homomorphisms defined by
\[
\<X^*\wedge Y^* \otimes P,\, Z^* \otimes Q\>_{\mu_v,\lambda_v,m} = s(X^*\wedge Y^*,Z^*) \cdot \<P,\iota_{m,v}(Q)\>_{\mu_v,\C}
\]
and
\[
\<W_1 \otimes X^*\wedge Y^* \otimes P,\, W_2 \otimes Z^* \otimes Q\>_{m,v} = Z_v(m+\tfrac{1}{2},W_1,W_2) \cdot s(X^*\wedge Y^*,Z^*) \cdot \<P,\iota_{m,v}(Q)\>_{\mu_v,\C}.
\]
Assume $\ell_v>\kappa_v$ for all $v \in S_\infty$ and $m+\tfrac{1}{2}$ is critical for $L(s,\itSigma \times \itPi)$. For $\sigma \in {\rm Aut}(\C)$, we then have the bilinear homomorphism
\[
\<\cdot,\cdot\>_{\sigma,m,\infty} = \bigotimes_{v \in S_\infty}\<\cdot,\cdot\>_{m,\sigma^{-1}\circ v}
\]
on 
\[
\left(\mathcal{W}({}^\sigma\!\itSigma_\infty,\psi_{3,\F_\infty}) \otimes \bigwedge^2 (\frak{g}_{3,\infty,\C} / \frak{k}_{3,\infty,\C})^* \otimes M_{{}^\sigma\!\mu,\C}^\vee\right)\times\left(\mathcal{W}({}^\sigma\!\itPi_\infty,\psi_{2,\F_\infty}) \otimes  (\frak{g}_{2,\infty,\C} / \frak{k}_{2,\infty,\C})^* \otimes M_{{}^\sigma\!\lambda,\C}^\vee\right).
\]
We simply write $\<\cdot,\cdot\>_{\sigma,m,\infty} = \<\cdot,\cdot\>_{m,\infty}$ when $\sigma$ is the identity map.

We have the following result due to Raghuram \cite{Raghuram2016} on the algebraicity of the critical values of $L(s,\itSigma \times \itPi)$ in terms of the bottom degree Whittaker periods.
Note that we have incorporated the period relation proved by Raghuram--Shahidi \cite{RS2008} into the formula. 

\begin{thm}[Raghuram]\label{T:Raghuram}
Assume $\ell_v > \kappa_v$ for all $v \in S_\infty$. Let $m+\tfrac{1}{2}$ be critical for $L(s,\itSigma \times \itPi)$ and $\underline{\varepsilon} \in \{\pm1\}^{S_\infty}$. We have
\begin{align*}
&\sigma\left(\frac{L^{(\infty)}(m+\tfrac{1}{2},\itSigma \times \itPi)}{|D_\F|^{1/2}\cdot G(\omega_\itPi)\cdot p^b(\itSigma)\cdot p(\itPi,\pm)}\cdot \<[\itSigma_\infty]_b , [\itPi_\infty]^{\underline{\varepsilon}}\>_{m,\infty}\right)\\
& = \frac{L^{(\infty)}(m+\tfrac{1}{2},{}^\sigma\!\itSigma \times {}^\sigma\!\itPi)}{|D_\F|^{1/2}\cdot G({}^\sigma\!\omega_\itPi)\cdot p^b({}^\sigma\!\itSigma)\cdot p({}^\sigma\!\itPi,\pm)}\cdot \<[{}^\sigma\!\itSigma_\infty]_b , [{}^\sigma\!\itPi_\infty]^{{}^\sigma\!\underline{\varepsilon}}\>_{\sigma,m,\infty}
\end{align*}
for all $\sigma \in {\rm Aut}(\C)$.
Here $[\itSigma_\infty]_b$ and $[\itPi_\infty]^{\underline{\varepsilon}}$ are the cohomology classes defined in (\ref{E:GL_3 class}) and (\ref{E:GL_2 class}), respectively.
\end{thm}

\begin{rmk}
The factor $|D_\F|^{1/2}$ is due to the comparison between Haar measures on $\A_\F$ that
\[
dx = |D_\F|^{1/2}\cdot dx^{\rm Tam},
\]
where $dx = \prod_v dx_v$ is the standard measure in \S\,\ref{SS:measure} and $dx^{\rm Tam}$ is the Tamagawa measure.
\end{rmk}

The non-vanishing of the archimedean factor $\<[\itSigma_\infty]_b , [\itPi_\infty]^{\underline{\varepsilon}}\>_{m,\infty}$ is determined by Kasten--Schmidt \cite{KS2013} and is also a special case of the result of Sun \cite{Sun2017} for $\GL_n \times \GL_{n-1}$. See also \cite[Theorem 3.8]{Schmidt1993} for the special case when $\ell=3$ and $\kappa=2$.
\begin{thm}[Kasten--Schmidt, Sun]\label{T:KS, Sun}
Assume $\ell_v > \kappa_v$ for all $v \in S_\infty$. Let $m+\tfrac{1}{2}$ be critical for $L(s,\itSigma \times \itPi)$ and $\underline{\varepsilon} = (\varepsilon_v)_{v \in S_\infty} \in \{\pm1\}^{S_\infty}$.
Then $\<[\itSigma_\infty]_b , [\itPi_\infty]^{\underline{\varepsilon}}\>_m \neq 0$ if and only if $\varepsilon_v = (-1)^m\epsilon_v(-1)$ for all $v \in S_\infty$.
\end{thm}

In Theorem \ref{T:RS} below, we refine Theorem \ref{T:Raghuram} by determine the rationality of the archimedean factor up to an explicit power of $2\pi\sqrt{-1}$.
Firstly we give a simplified formula for the archimedean factor in the following lemma.
\begin{lemma}\label{L:archimedean factor 1}
Let $v \in S_\infty$. Assume $\ell_v > \kappa_v$ and $m+\tfrac{1}{2}$ is is not a pole of neither $L(s,\itSigma_v \times \itPi_v)$ nor $L(1-s,\itSigma_v^\vee \times \itPi_v^\vee)$.
We have
\begin{align*}
& \<[\itSigma_v]_b,[\itPi_v]^\pm\>_{m,v}\\
& = \frac{(\sqrt{-1})^{{\sf w}(\itSigma)/2}}{(\ell_v-\kappa_v)!}\cdot Z_v(m+\tfrac{1}{2},W_{(\ell_v, {\sf w}(\itSigma), \epsilon_v;\,-\kappa_v)}, W_{\kappa_v}^+ \otimes |\mbox{ }|^{{\sf w}(\itPi)/2})\\
& \times\left\<E_-^{\ell_v-\kappa_v}\cdot \left(X_{-1}^*\wedge X_{-2}^* \otimes \rho_{\mu_v^\vee}\left( \bp 1 & 0 & 1  \\ \sqrt{-1} & 0 & -\sqrt{-1} \\ 0 & 1 & 0\ep\right)P_{\mu_v^\vee}^+\right), \,Y_+^* \otimes (\sqrt{-1}\,x+y)^{\kappa_v-2}\right\>_{\mu_v,\lambda_v,m}\\
& \pm (-1)^{\kappa_v+{\sf w}(\itSigma)/2}\frac{(\sqrt{-1})^{{\sf w}(\itSigma)/2+{\sf w}(\itPi)}}{(\ell_v-\kappa_v)!}\cdot Z_v(m+\tfrac{1}{2},W_{(\ell_v, {\sf w}(\itSigma), \epsilon_v;\,\kappa_v)}, W_{\kappa_v}^- \otimes |\mbox{ }|^{{\sf w}(\itPi)/2})\\
& \times \left\<E_+^{\ell_v-\kappa_v}\cdot \left(X_{1}^*\wedge X_{2}^* \otimes \rho_{\mu_v^\vee}\left( \bp -1 & 0 & -1  \\ \sqrt{-1} & 0 & -\sqrt{-1} \\ 0 & 1 & 0\ep\right)P_{\mu_v^\vee}^+\right), Y_-^* \otimes (x+\sqrt{-1}\,y)^{\kappa_v-2}\right\>_{\mu_v,\lambda_v,m}.
\end{align*}
\end{lemma}

\begin{proof}
It is clear that the local zeta integral defines an element in 
${\rm Hom}_{\GL_2(\R)}(\itSigma_v \otimes \itPi_v, \,|\mbox{ }|^{-s+1/2})$.
Hence
$Z_v(m+\tfrac{1}{2},W_{(\ell_v,{\sf w}(\itSigma),\epsilon_v;\,i)},W_{\kappa_v}^\pm\otimes|\mbox{ }|^{{\sf w}(\itPi)/2}) \neq 0
$ only when $i=\mp \kappa_v$.
The assertion then follows from (\ref{E:relation}) of Lemma \ref{L:cohomology generator} and the definition of $\<\cdot,\cdot\>_m$.
\end{proof}

In the following lemma, we compute the archimedean local zeta integrals appearing in Lemma \ref{L:archimedean factor 1}.
The result was announced in \cite[Theorem 4.1]{HIM2016} and we give a proof for the sake of completeness.
\begin{lemma}\label{L:archimedean RS}
Let $v \in S_\infty$. Assume $\ell_v \geq \kappa_v$.
We have
\begin{align*}
Z_v(s,W_{(\ell_v, {\sf w}(\itSigma), \epsilon_v;\,\kappa_v)},W_{\kappa_v}^-\otimes |\mbox{ }|^{{\sf w}(\itPi)/2}) &= (-1)^{{\sf w}(\itSigma)/2}\epsilon_v(-1)\cdot Z_v(s,W_{(\ell_v, {\sf w}(\itSigma), \epsilon_v;\,-\kappa_v)},W_{\kappa_v}^+\otimes |\mbox{ }|^{{\sf w}(\itPi)/2})\\
& = (\sqrt{-1})^{2\kappa_v-\ell_v} \cdot L(s,\itSigma_v \times \itPi_v).
\end{align*}
\end{lemma}

\begin{proof}
We drop the subscript $v$ for brevity.
We may assume ${\sf w}(\itSigma)={\sf w}(\itPi)=0$ after replacing $\itSigma_v$ and $\itPi_v$ by $\itSigma_v \otimes |\mbox{ }|^{-{\sf w}(\itSigma)/2}$ and $\itPi_v \otimes |\mbox{ }|^{-{\sf w}(\itPi)/2}$, respectively. Write $W_{(\ell;\, \pm\kappa)} = W_{(\ell, {\sf w}(\itSigma), \epsilon;\,\pm\kappa)}$ and $W_{(\ell;\, (j_1,j_2,j_3))} = W_{(\ell, {\sf w}(\itSigma), \epsilon;\,(j_1,j_2,j_3))}$.
Since the assertion is an equality for meromorphic functions in $s \in \C$, it suffices to prove the equality for ${\rm Re}(s)$ sufficiently large.
We assume ${\rm Re}(s)>\tfrac{\kappa+1}{2}$. Let $L_1$ and $L_2$ be vertical paths from $c_1-\sqrt{-1}\,\infty$ to $c_1+\sqrt{-1}\,\infty$ and $1-\sqrt{-1}\,\infty$ to $1+\sqrt{-1}\,\infty$, respectively, where 
$\kappa<c_1<{\rm Re}(s)+\tfrac{\kappa-1}{2}$.
We have
\begin{align*}
&Z(s,W_{(\ell;\, \kappa)},W_\kappa^-)\\
& = \int_{N_2(\R)\backslash\GL_2(\R)} W_{(\ell;\, \kappa)}(\iota(g))W_\kappa^-\left({\rm diag}(-1,1) g\right)|\det(g)|^{s-1/2}\,dg\\
& = \int_{0}^\infty d^\times a_1 \int_{0}^\infty d^\times a_2\, a_1^{s-3/2}a_2^{2s-1}W_{(\ell;\, \kappa)}\left(\iota\left({\rm diag}(a_1a_2,a_2)\right)\right)W_\kappa^-\left({\rm diag}(-a_1a_2,a_2)\right)\\
& = \sum_{i=0}^{\kappa}(\sqrt{-1})^{\kappa-i}{\kappa \choose i}\int_0^\infty \int_{0}^\infty a_1^{s+(\kappa-1)/2}a_2^{2s}e^{-2\pi a_1}W_{(\ell;\, (i,\kappa-i,\ell-\kappa))}\left(\iota\left({\rm diag}(a_1a_2,a_2)\right)\right)\,d^\times a_1 d^\times a_2\\
& = (\sqrt{-1})^{2\kappa-\ell}
\sum_{i=0}^{\kappa}{\kappa \choose i}\int_0^\infty d^\times a_1\int_{0}^\infty d^\times a_2\int_{L_1}\frac{ds_1}{2\pi\sqrt{-1}}\int_{L_2}\frac{ds_2}{2\pi\sqrt{-1}}\, a_1^{-s_1+s+(\kappa-1)/2}a_2^{-s_2+2s}e^{-2\pi a_1}\\
&\quad\quad \quad\quad\quad\quad\quad\quad\quad\quad\quad\quad\quad\quad\quad\quad\quad\quad\times\frac{\Gamma_\C(s_1+\tfrac{\ell-1}{2})\Gamma_\R(s_1+i)\Gamma_\C(s_2+\tfrac{\ell-1}{2})\Gamma_\R(s_2+\ell-\kappa)}{\Gamma_\R(s_1+s_2+\ell-\kappa+i)}\\
& = 2^{-1}(\sqrt{-1})^{2\kappa-\ell}
\sum_{i=0}^{\kappa}{\kappa \choose i}\int_0^\infty d^\times a_2\int_{L_1}\frac{ds_1}{2\pi\sqrt{-1}}\int_{L_2}\frac{ds_2}{2\pi\sqrt{-1}}\,a_2^{-s_2+2s}\\&\quad\quad\quad\quad\quad\quad\quad\quad\quad\quad\quad\quad\times\frac{\Gamma_\C(s_1+\tfrac{\ell-1}{2})\Gamma_\R(s_1+i)\Gamma_\C(-s_1+s+\tfrac{\kappa-1}{2})\Gamma_\C(s_2+\tfrac{\ell-1}{2})\Gamma_\R(s_2+\ell-\kappa)}{\Gamma_\R(s_1+s_2+\ell-\kappa+i)}\\
& = 2^{-1}(\sqrt{-1})^{2\kappa-\ell}
\Gamma_\C(2s+\tfrac{\ell-1}{2})\Gamma_\R(2s+\ell-\kappa)\\
&\times\sum_{i=0}^{\kappa}{\kappa \choose i} \int_{L_1}\frac{ds_1}{2\pi\sqrt{-1}}\,\frac{\Gamma_\C(s_1+\tfrac{\ell-1}{2})\Gamma_\R(s_1+i)\Gamma_\C(-s_1+s+\tfrac{\kappa-1}{2})}{\Gamma_\R(s_1+2s+\ell-\kappa+i)}.
\end{align*}
Here the last equality follows from the Mellin inversion formula.
Now we make a change of variable from $s_1$ to $s_1-i$.
By \cite[Lemma 12.7-(i)]{Ichino2005}, we have
\[
\sum_{i=0}^{\kappa}{\kappa \choose i}\Gamma_\C(s_1+\tfrac{\ell-1}{2}-i)\Gamma_\C(-s_1+s+\tfrac{\kappa-1}{2}+i) = \frac{\Gamma_\C(-s_1+s+\tfrac{\kappa-1}{2})\Gamma_\C(s+\tfrac{\ell+\kappa}{2}-1)\Gamma_\C(s_1+\tfrac{\ell-1}{2}-\kappa)}{\Gamma_\C(s+\tfrac{\ell-\kappa}{2}-1)}.
\]
Therefore, by our choice of the path $L_1$, we have
\begin{align*}
Z(s,W_{(\ell;\, \kappa)},W_\kappa^-) &= 2^{-1}(\sqrt{-1})^{2\kappa-\ell}\cdot\frac{\Gamma_\C(2s+\tfrac{\ell-1}{2})\Gamma_\R(2s+\ell-\kappa)\Gamma_\C(s+\tfrac{\ell+\kappa}{2}-1)}{\Gamma_\C(s+\tfrac{\ell-\kappa}{2}-1)}\\
&\times\int_{L_1}\frac{ds_1}{2\pi\sqrt{-1}}\,\frac{\Gamma_\C(s_1+\tfrac{\ell-1}{2}-\kappa)\Gamma_\R(s_1)\Gamma_\C(-s_1+s+\tfrac{\kappa-1}{2})}{\Gamma_\R(s_1+2s+\ell-\kappa)}.
\end{align*}
By the second Barnes lemma, we have
\begin{align*}
&\int_{L_1}\frac{ds_1}{2\pi\sqrt{-1}}\,\frac{\Gamma_\C(s_1+\tfrac{\ell-1}{2}-\kappa)\Gamma_\R(s_1)\Gamma_\C(-s_1+s+\tfrac{\kappa-1}{2})}{\Gamma_\R(s_1+2s+\ell-\kappa)}\\
& = 2 \cdot\frac{\Gamma_\C(s+\tfrac{\ell-\kappa}{2}-1)\Gamma_\C(s+\tfrac{\ell-\kappa}{2})\Gamma_\C(s+\tfrac{\kappa-1}{2})}{\Gamma_\C(2s+\tfrac{\ell-1}{2})\Gamma_\R(2s+\ell-\kappa)}.
\end{align*}
We conclude that
\begin{align*}
Z(s,W_{(\ell;\, \kappa)},W_\kappa^-) &= (\sqrt{-1})^{2\kappa-\ell}\cdot 
\Gamma_\C(s+\tfrac{\ell+\kappa}{2}-1)\Gamma_\C(s+\tfrac{\ell-\kappa}{2})\Gamma_\C(s+\tfrac{\kappa-1}{2})\\
& = (\sqrt{-1})^{2\kappa-\ell}\cdot 
 L(s,\itSigma_v \times \itPi_v).
\end{align*}
Note that $W_{(\ell;\, -\kappa)}$ and $W_\kappa^+$ are equal to the right translations of $(-1)^{{\sf w}(\itSigma)/2}\epsilon(-1)\cdot W_{(\ell;\, \kappa)}$ and $W_\kappa^-$ by ${\rm diag}(-1,1,1)$ and ${\rm diag}(-1,1)$, respectively.
Hence
\[
Z(s,W_{(\ell;\, \kappa)},W_\kappa^-) = (-1)^{{\sf w}(\itSigma)/2}\epsilon(-1)\cdot Z(s,W_{(\ell;\, -\kappa)},W_\kappa^+).
\]
This completes the proof.
\end{proof}

The following lemma will be used in the proof of Theorem \ref{T:RS}.
\begin{lemma}\label{L:combinatorial}
Let $v \in S_\infty$.
The Lie algebra action of ${\rm Ad}\left(\bp 1 & 1  & 0\\ \sqrt{-1} & -\sqrt{-1} & 0 \\ 0 & 0 & 1\ep\right)^{-1}E_\pm$ on $M_{\mu_v,\C}^\vee$ is defined over $\Q$.
\end{lemma}

\begin{proof}
Note that
\[
{\rm Ad}\left(\bp 1 & 1  & 0\\ \sqrt{-1} & -\sqrt{-1} & 0 \\ 0 & 0 & 1\ep\right)^{-1}E_\pm = \frac{1}{2}\bp 0&0&1 \\ 0&0&1 \\ -1&-1&0\ep \pm \frac{1}{2}\bp 0&0&-\sqrt{-1} \\ 0&0&\sqrt{-1} \\ -\sqrt{-1}&\sqrt{-1}&0\ep \otimes \sqrt{-1}.
\]
It is clear that $\bp 0&0&1 \\ 0&0&1 \\ -1&-1&0\ep\cdot M_{\mu_v}^\vee \subset M_{\mu_v}^\vee$.
Since
\begin{align*}
&\bp x_{11} & x_{12} & x_{13} \\ x_{21} & x_{22} & x_{23}\ep \exp\left( t \bp 0&0&0 \\ 0&0&\sqrt{-1} \\ 0&\sqrt{-1}&0\ep \right) \\
&= \bp x_{11} & x_{12}\cos t +\sqrt{-1}\,x_{13}\sin t&  \sqrt{-1}\,x_{12}\sin t + x_{13}\cos t\\ x_{21} & x_{22}\cos t +\sqrt{-1}\,x_{23}\sin t&  \sqrt{-1}\,x_{22}\sin t + x_{23}\cos t\ep,
\end{align*}
we see that $\bp 0&0&0 \\ 0&0&\sqrt{-1} \\ 0&\sqrt{-1}&0\ep\cdot M_{\mu_v}^\vee \subset \sqrt{-1}\cdot M_{\mu_v}^\vee$.
Similarly $\bp 0&0&\sqrt{-1} \\ 0&0&0 \\ \sqrt{-1}&0&0\ep\cdot M_{\mu_v}^\vee \subset \sqrt{-1}\cdot M_{\mu_v}^\vee$.
This completes the proof.
\end{proof}

\begin{thm}\label{T:RS}
Assume $\ell_v > \kappa_v$ for all $v \in S_\infty$. Let $m+\tfrac{1}{2}$ be critical for $L(s,\itSigma \times \itPi)$.
We have
\begin{align*}
&\sigma\left(\frac{L^{(\infty)}(m+\tfrac{1}{2},\itSigma \times \itPi)}{|D_\F|^{1/2}\cdot(2\pi\sqrt{-1})^{3dm+\sum_{v \in S_\infty}(2\ell_v+\kappa_v+3{\sf w}(\itSigma)+3{\sf w}(\itPi))/2}\cdot G(\omega_\itPi)\cdot p^b(\itSigma)\cdot p(\itPi,(-1)^m\cdot\underline{\varepsilon}(\itSigma))}\right) \\
&= \frac{L^{(\infty)}(m+\tfrac{1}{2},{}^\sigma\!\itSigma \times {}^\sigma\!\itPi)}{|D_\F|^{1/2}\cdot(2\pi\sqrt{-1})^{3dm+\sum_{v \in S_\infty}(2\ell_v+\kappa_v+3{\sf w}(\itSigma)+3{\sf w}(\itPi))/2}\cdot G({}^\sigma\!\omega_\itPi)\cdot p^b({}^\sigma\!\itSigma)\cdot p({}^\sigma\!\itPi,(-1)^m\cdot\underline{\varepsilon}(\itSigma))}
\end{align*}
for all $\sigma \in {\rm Aut}(\C)$.
Here $\underline{\varepsilon}(\itSigma) = (\epsilon_v(-1))_{v\in S_\infty}$.
\end{thm}

\begin{proof}
By the definition of the classes $[\itSigma_\infty]_b$, $[{}^\sigma\!\itSigma_\infty]_b$, $[\itPi_\infty]^{\underline{\varepsilon}}$, and $[{}^\sigma\!\itPi_\infty]^{{}^\sigma\!\underline{\varepsilon}}$, we have
\[
\<[{}^\sigma\!\itSigma_\infty]_b , [{}^\sigma\!\itPi_\infty]_b^{{}^\sigma\!\underline{\varepsilon}}\>_{\sigma,m,\infty} = \<[\itSigma_\infty]_b , [\itPi_\infty]_b^{\underline{\varepsilon}}\>_{m,\infty}
\]
for all $\underline{\varepsilon} \in \{\pm 1\}^{S_\infty}$.
Indeed, 
\begin{align*}
\<[{}^\sigma\!\itSigma_\infty]_b , [{}^\sigma\!\itPi_\infty]_b^{{}^\sigma\!\underline{\varepsilon}}\>_{\sigma,m,\infty}& = \prod_{v \in S_\infty}\<[{}^\sigma\!\itSigma_v]_b , [{}^\sigma\!\itPi_v]_b^{\varepsilon_{\sigma^{-1}\circ v}}\>_{m,\sigma^{-1}\circ v}\\
& =  \prod_{v \in S_\infty}\<[{}^\sigma\!\itSigma_{\sigma\circ v}]_b , [{}^\sigma\!\itPi_{\sigma\circ v}]_b^{\varepsilon_v}\>_{m,v}\\
& = \prod_{v \in S_\infty}\<[\itSigma_v]_b , [\itPi_v]_b^{\varepsilon_v}\>_{m,v}\\
& = \<[\itSigma_\infty]_b , [\itPi_\infty]_b^{\underline{\varepsilon}}\>_{m,\infty}.
\end{align*}
By Theorems \ref{T:Raghuram} and \ref{T:KS, Sun}, it suffices to show that
\[
\<[\itSigma_v]_b , [\itPi_v]^{\pm}\>_m \in (2\pi\sqrt{-1})^{-3m-(2\ell_v+\kappa_v+3{\sf w}(\itSigma)+3{\sf w}(\itPi))/2} \cdot \Q
\]
for all $v \in S_\infty$.
Fix $v \in S_\infty$. We drop the subscript $v$ for brevity.
Note that 
\begin{align*}
(\sqrt{-1}\,x+y)^{\kappa-2} &= (-2)^{m+(\kappa+{\sf w}(\itPi))/2}(\sqrt{-1})^{-2+m+(3\kappa+{\sf w}(\itPi))/2}\cdot\rho_{\lambda^\vee-m}\left( \bp 1 & 1 \\ \sqrt{-1} & -\sqrt{-1}\ep\right) y^{\kappa-2},\\
(\sqrt{-1})^{{\sf w}(\itPi)}(x+\sqrt{-1}\,y)^{\kappa-2} &= (-2)^{m+(\kappa+{\sf w}(\itPi))/2}(\sqrt{-1})^{m+(\kappa+3{\sf w}(\itPi))/2}\cdot\rho_{\lambda^\vee-m}\left( \bp 1 & 1 \\ \sqrt{-1} & -\sqrt{-1}\ep\right) x^{\kappa-2}.
\end{align*}
It is easy to see that the kernels of the linear functionals $s(\cdot,Y_+^*)$ and $s(\cdot,Y_-^*)$ do not contain $X_0^* \wedge X_{-2}^*$ and $X_0^* \wedge X_2^*$, respectively. 
Indeed, we have $s(X_0^*\wedge X_{-2}^*,Y_+^*) = 8$ and $s(X_0^*\wedge X_2^*,Y_-^*) = -8$.
Thus we only need to consider the $X_0^* \wedge X_{-2}^*$ and $X_0^* \wedge X_2^*$ components of 
\[
E_-^{\ell-\kappa}\cdot \left(X_{-1}^*\wedge X_{-2}^* \otimes \rho_{\mu^\vee}\left( \bp 1 & 0 & 1  \\ \sqrt{-1} & 0 & -\sqrt{-1} \\ 0 & 1 & 0\ep\right)P_{\mu^\vee}^+\right)
\]
and
\[
E_+^{\ell-\kappa}\cdot \left(X_{1}^*\wedge X_{2}^* \otimes \rho_{\mu^\vee}\left( \bp -1 & 0 & -1  \\ \sqrt{-1} & 0 & -\sqrt{-1} \\ 0 & 1 & 0\ep\right)P_{\mu^\vee}^+\right),
\]
respectively.
By (\ref{E:Lie algebra action 2}), these components are equal to 
\[
2\cdot X_0^* \wedge X_{-2}^* \otimes E_-^{\ell-\kappa-1}\cdot \rho_{\mu^\vee}\left( \bp 1 & 0 & 1  \\ \sqrt{-1} & 0 & -\sqrt{-1} \\ 0 & 1 & 0\ep\right)P_{\mu^\vee}^+
\]
and
\[
-2\cdot X_0^* \wedge X_{2}^* \otimes E_+^{\ell-\kappa-1}\cdot \rho_{\mu^\vee}\left( \bp -1 & 0 & -1  \\ \sqrt{-1} & 0 & -\sqrt{-1} \\ 0 & 1 & 0\ep\right)P_{\mu^\vee}^+,
\]
respectively.
We conclude that
\begin{align*}
& \left\<E_-^{\ell-\kappa}\cdot \left(X_{-1}^*\wedge X_{-2}^* \otimes \rho_{\mu^\vee}\left( \bp 1 & 0 & 1  \\ \sqrt{-1} & 0 & -\sqrt{-1} \\ 0 & 1 & 0\ep\right)P_{\mu^\vee}^+\right), Y_+^* \otimes (\sqrt{-1}\,x+y)^{\kappa-2}\right\>_{\mu,\lambda,m}\\
& = (-2)^{4+m+(\kappa+{\sf w}(\itPi))/2}(\sqrt{-1})^{-2+m+(3\kappa+{\sf w}(\itPi))/2} \\
&\times \left\<E_-^{\ell-\kappa-1}\cdot \rho_{\mu^\vee}\left( \bp 1 & 0 & 1  \\ \sqrt{-1} & 0 & -\sqrt{-1} \\ 0 & 1 & 0\ep\right)P_{\mu^\vee}^+ , \,\rho_{\mu}\left( \bp 1 & 1 &0\\ \sqrt{-1} & -\sqrt{-1} & 0 \\ 0 & 0 & 1\ep\right) \iota_m(y^{\kappa-2})\right\>_{\mu,\C}\\
& = (-2)^{4+m+(\kappa+{\sf w}(\itPi))/2}(\sqrt{-1})^{-2+m+(3\kappa+{\sf w}(\itPi))/2}\\
&\times\left\<{\rm Ad}\left(\bp 1 & 1  & 0\\ \sqrt{-1} & -\sqrt{-1} & 0 \\ 0 & 0 & 1\ep\right)^{-1}E_-^{\ell-\kappa-1}\cdot \rho_{\mu^\vee}\left( \bp 1 & 0 & 0  \\ 0 & 0 & 1 \\ 0 & 1 & 0\ep\right)P_{\mu^\vee}^+ , \,\iota_m(y^{\kappa-2})\right\>_{\mu,\C}\\
& \in (\sqrt{-1})^{m+(3\kappa+{\sf w}(\itPi))/2}\cdot\Q.
\end{align*}
Here the last line follows from Lemma \ref{L:combinatorial}.
Similarly we have
\begin{align*}
& (\sqrt{-1})^{{\sf w}(\itPi)}\left\<E_+^{\ell-\kappa}\cdot \left(X_1^*\wedge X_2^* \otimes \rho_{\mu^\vee}\left( \bp -1 & 0 & -1  \\ \sqrt{-1} & 0 & -\sqrt{-1} \\ 0 & 1 & 0\ep\right)P_{\mu^\vee}^+\right), Y_-^* \otimes (x+\sqrt{-1}\,y)^{\kappa-2}\right\>_{\mu,\lambda,m}\\
& = (-2)^{4+m+(\kappa+{\sf w}(\itPi))/2}(\sqrt{-1})^{m+(\kappa+3{\sf w}(\itPi))/2}\\
&\times\left\<{\rm Ad}\left(\bp 1 & 1  & 0\\ \sqrt{-1} & -\sqrt{-1} & 0 \\ 0 & 0 & 1\ep\right)^{-1}E_+^{\ell-\kappa-1}\cdot \rho_{\mu^\vee}\left( \bp 0 & 0 & -1  \\ -1 & 0 & 0 \\ 0 & 1 & 0\ep\right)P_{\mu^\vee}^+ , \,\iota_m(x^{\kappa-2})\right\>_{\mu,\C}\\
& \in (\sqrt{-1})^{m+(\kappa+3{\sf w}(\itPi))/2}\cdot\Q.
\end{align*}
Finally, by Lemma \ref{L:archimedean RS}, we have
\begin{align*}
&(\sqrt{-1})^{{\sf w}(\itSigma)/2}Z(m+\tfrac{1}{2},W_{(\ell, {\sf w}(\itSigma), \epsilon;\,\kappa)},W_\kappa^-\otimes |\mbox{ }|^{{\sf w}(\itPi)/2}) \\
&= (-\sqrt{-1})^{{\sf w}(\itSigma)/2}\epsilon(-1)\cdot Z(m+\tfrac{1}{2},W_{(\ell, {\sf w}(\itSigma), \epsilon;\,-\kappa)},W_\kappa^+\otimes |\mbox{ }|^{{\sf w}(\itPi)/2})\\
& \in (\sqrt{-1})^{-\ell+{\sf w}(\itSigma)/2}\cdot \pi^{-3m-(2\ell+\kappa+3{\sf w}(\itSigma)+3{\sf w}(\itPi))/2} \cdot \Q^\times.
\end{align*}
This completes the proof.
\end{proof}

\begin{rmk}
The theorem is compatible with the result of Januszewski \cite{Januszewski2019} for $\GL_n \times \GL_{n-1}$.
\end{rmk}

\subsection{A period relation for the Gelbart--Jacquet lifting}

We consider the special case when $\itSigma = {\rm Sym}^2\itPi$ is the Gelbart--Jacquet lift of $\itPi$ (cf.\,\cite{GJ1978}), which is cuspidal by our assumption that $\itPi$ is non-CM. Then ${\rm Sym}^2\itPi$ is cohomological with
\[
\ell_v=2\kappa_v-1,\quad {\sf w}(\itSigma) = 2{\sf w}(\itPi),\quad \epsilon_v(-1) = (-1)^{1+{\sf w}(\itPi)}
\]
for all $v \in S_\infty$.
In the following theorem, we prove a period relation between the bottom degree Whittaker period $p^b({\rm Sym}^2\itPi)$ and the Petersson norm $\Vert f_\itPi \Vert$.

\begin{thm}\label{T:GJ period relation}
Suppose $\kappa_v \geq 3$ for all $v \in S_\infty$.
Then we have
\[
\sigma\left(\frac{|D_\F|^{1/2}\cdot \pi^{2d}\cdot G(\omega_\itPi)^3\cdot\Vert f_\itPi \Vert^2}{p^b({\rm Sym}^2\itPi)}\right) = \frac{|D_\F|^{1/2}\cdot \pi^{2d}\cdot G({}^\sigma\!\omega_\itPi)^3\cdot\Vert f_{{}^\sigma\!\itPi} \Vert^2}{p^b({\rm Sym}^2 {}^\sigma\!\itPi)}
\]
for all $\sigma \in {\rm Aut}(\C)$.
\end{thm}

\begin{proof}

Let $L(s,\itPi \times \itPi \times \itPi)$ be the triple product $L$-function of $\itPi \times \itPi \times \itPi$ (cf.\,\cite{Garrett1987}, \cite{PSR1987}, and \cite{Ikeda1989}).
Let $m+\tfrac{1}{2} \in \Z+\tfrac{1}{2}$ be a critical point for $L(s,\itPi \times \itPi \times \itPi)$.
By the result of Garrett--Harris \cite{GH1993} and Theorem \ref{T:RT2011}, we have
\begin{align*}
&\sigma\left(\frac{L^{(\infty)}(m+\tfrac{1}{2},\itPi \times \itPi \times \itPi)}{(2\pi\sqrt{-1})^{4dm}\cdot \pi^{2d+3\sum_{v \in S_\infty}\kappa_v+6d{\sf w}(\itPi)}\cdot G(\omega_\itPi)^6 \cdot \Vert f_\itPi \Vert^3}\right) \\
&= \frac{L^{(\infty)}(m+\tfrac{1}{2},{}^\sigma\!\itPi \times {}^\sigma\!\itPi \times {}^\sigma\!\itPi)}{(2\pi\sqrt{-1})^{4dm}\cdot \pi^{2d+3\sum_{v \in S_\infty}\kappa_v+6{\sf w}(\itPi)}\cdot G({}^\sigma\!\omega_\itPi)^6 \cdot \Vert f_{{}^\sigma\!\itPi} \Vert^3}
\end{align*}
and 
\begin{align*}
&\sigma\left( \frac{L^{(\infty)}(m+\tfrac{1}{2},\itPi \times \omega_\itPi)}{(2\pi\sqrt{-1})^{dm+\sum_{v \in S_\infty}(\kappa_v+3{\sf w}(\itPi))/2}\cdot G(\omega_\itPi) \cdot p(\itPi,(-1)^m\cdot {\rm sgn}(\omega_\itPi))}\right)\\
& = \frac{L^{(\infty)}(m+\tfrac{1}{2},{}^\sigma\!\itPi \times {}^\sigma\!\omega_\itPi)}{(2\pi\sqrt{-1})^{dm+\sum_{v \in S_\infty}(\kappa_v+3{\sf w}(\itPi))/2}\cdot G({}^\sigma\!\omega_\itPi) \cdot p({}^\sigma\!\itPi,(-1)^m\cdot {\rm sgn}(\omega_\itPi))},
\end{align*}
respectively, for all $\sigma \in {\rm Aut}(\C)$.
On the other hand, we have the factorization of $L$-functions
\[
L(s,\itPi \times \itPi \times \itPi) = L(s,{\rm Sym}^2\itPi \times \itPi)\cdot L(s,\itPi \times \omega_\itPi).
\]
Now we take 
\[
m+\tfrac{1}{2} = \tfrac{\min_{v \in S_\infty}\{\kappa_v\}-3{\sf w}(\itPi)-1}{2}
\]
be the rightmost critical point. Then the assumption $\kappa \geq 3$ implies that the above $L$-functions are all non-vanishing at $m+\tfrac{1}{2}$.
Combining with Lemma \ref{L:period relation GL_2}, we conclude that
\begin{align*}
&\sigma \left( \frac{L^{(\infty)}(m+\tfrac{1}{2},{\rm Sym}^2\itPi \times \itPi)}{(2\pi\sqrt{-1})^{3dm+\sum_{v \in S_\infty}(5\kappa_v+9{\sf w}(\itPi)+2)/2}\cdot G(\omega_\itPi)^4\cdot \Vert f_\itPi \Vert^2 \cdot p(\itPi,(-1)^{m+1}\cdot {\rm sgn}(\omega_\itPi))}\right)\\
& = \frac{L^{(\infty)}(m+\tfrac{1}{2},{\rm Sym}^2{}^\sigma\!\itPi \times {}^\sigma\!\itPi)}{(2\pi\sqrt{-1})^{3dm+\sum_{v \in S_\infty}(5\kappa_v+9{\sf w}(\itPi)+2)/2}\cdot G({}^\sigma\!\omega_\itPi)^4\cdot \Vert f_{{}^\sigma\!\itPi} \Vert^2 \cdot p({}^\sigma\!\itPi,(-1)^{m+1}\cdot {\rm sgn}(\omega_\itPi))}
\end{align*}
for all $\sigma \in {\rm Aut}(\C)$.
Note that
\[
-{\rm sgn}(\omega_\itPi) = ((-1)^{1+\kappa_v})_{v \in S_\infty} = (\epsilon_v(-1))_{v \in S_\infty}.
\]
Comparing with Theorem \ref{T:RS}, we then obtain the desired period relation.
This completes the proof.
\end{proof}

\subsection{A period relation for $\GL_3$}\label{SS:Deligne}

We fix the notation as in the beginning of \S\,\ref{S:GL_3 GL_2}. We assume $\F=\Q$ in this section for simplicity.
Let $\mathbb{M}_\itSigma$ and $\mathbb{M}_\itPi$ be the (conjectural) motives attached to $\itSigma$ and $\itPi$, respectively, proposed by Clozel \cite[Conjecture 4.5]{Clozel1990}. 
Note that the existence of $\mathbb{M}_\itPi$ was proved by \cite{Scholl1990}. 
Then $\mathbb{M}_\itSigma$ and $\mathbb{M}_\itPi$ are pure motives over $\Q$ with coefficients in $\Q(\itSigma)$ and $\Q(\itPi)$, of rank $3$ and $2$, and weights $-{\sf w}(\itSigma)$ and $-{\sf w}(\itPi)-1$, respectively, and
\[
L(s,\mathbb{M}_\itSigma) = (L^{(\infty)}(s,{}^\sigma\!\itSigma))_\sigma,\quad L(s,\mathbb{M}_\itPi) = (L^{(\infty)}(s+\tfrac{1}{2},{}^\sigma\!\itPi))_\sigma.
\]
Here $\sigma$ runs through complete coset representatives of ${\rm Aut}(\C) / {\rm Aut}(\C / \Q(\itSigma))$ and ${\rm Aut}(\C) / {\rm Aut}(\C / \Q(\itPi))$, respectively.
Let 
\[
c^\pm(\mathbb{M}_\itSigma) \in (\Q(\itSigma)\otimes_\Q \C)^\times,\quad c^\pm(\mathbb{M}_\itPi) \in (\Q(\itPi)\otimes_\Q \C)^\times
\]
be Deligne's periods, well-defined up to $\Q(\itSigma)^\times$ and $\Q(\itPi)^\times$, for the motives $\mathbb{M}_\itSigma$ and $\mathbb{M}_\itPi$, respectively. 
Consider the tensor motive $\mathbb{M}_\itSigma \otimes \mathbb{M}_\itPi$, which is a pure motive over $\Q$ with coefficients in $\Q(\itSigma,\itPi) = \Q(\itSigma)\cdot\Q(\itPi)$, of rank $6$, weight $-{\sf w}(\itSigma)-{\sf w}(\itPi)-1$, and 
\[
L(s,\mathbb{M}_\itSigma \otimes \mathbb{M}_\itPi) = (L^{(\infty)}(s+\tfrac{1}{2}, {}^\sigma\!\itSigma \times {}^\sigma\!\itPi))_\sigma.
\]
In \cite{Yoshida2001}, Yoshida computed the Deligne's periods for tensor product motives in terms of fundamental periods. The result of Yoshida was explicated further by Bhagwat \cite{Bhagwat2014}. Specializing \cite[Theorem 3.2]{Bhagwat2014} to our tensor motive $\mathbb{M}_\itSigma \otimes \mathbb{M}_\itPi$, under the assumption that $\ell > \kappa$, we have
\[
c^\pm(\mathbb{M}_\itSigma \otimes \mathbb{M}_\itPi) = (2\pi\sqrt{-1})^{{\sf w}(\itPi)+1}\cdot(G({}^\sigma\!\omega_\itPi))_\sigma \cdot c^+({\mathbb M}_\itSigma)\cdot c^-({\mathbb M}_\itSigma) \cdot c^{\pm\epsilon(-1)}(\mathbb{M}_\itPi).
\]
Here we have enlarge the coefficients to $\Q(\itSigma,\itPi)$.
Deligne's conjecture then predicts that
\begin{align}\label{E:Deligne conj.}
\frac{L(m,\mathbb{M}_\itSigma \otimes \mathbb{M}_\itPi)}{(2\pi\sqrt{-1})^{3m+{\sf w}(\itPi)+1}\cdot(G({}^\sigma\!\omega_\itPi))_\sigma \cdot c^+({\mathbb M}_\itSigma)\cdot c^-({\mathbb M}_\itSigma) \cdot c^{(-1)^m\epsilon(-1)}(\mathbb{M}_\itPi)} \in \Q(\itSigma,\itPi)
\end{align}
for all critical points $m \in \Z$ of $\mathbb{M}_\itSigma \otimes \mathbb{M}_\itPi$.
On the other hand, by Theorem \ref{T:RT2011} we have
\[
c^\pm(\mathbb{M}_\itPi) = (2\pi\sqrt{-1})^{(\kappa+{\sf w}(\itPi))/2}\cdot (p({}^\sigma\!\itPi,\pm))_\sigma.
\]
Therefore, comparing (\ref{E:Deligne conj.}) with Theorem \ref{T:RS}, it is natural to propose the following conjecture.

\begin{conj}\label{C:Deligne-Whittaker}
We have
\[
c^+({\mathbb M}_\itSigma)\cdot c^-({\mathbb M}_\itSigma) = (2\pi\sqrt{-1})^{\ell+3{\sf w}(\itSigma)/2-1}\cdot\left(p^b({}^\sigma\!\itSigma)\right)_\sigma \in (\Q(\itSigma)\otimes_\Q \C)^\times.
\]
\end{conj}

\begin{rmk}
It would be interesting to compare the conjecture with the results \cite{Mahnkopf2005}, \cite{RS2017}, and \cite{Sachdeva2020} on the algebraicity of the twisted standard $L$-functions of $\itSigma$.
\end{rmk}

\section{Algebraicity of the adjoint $L$-functions for $\GL_3$}\label{S:adjoint}

Let $\itSigma = \bigotimes_v \itSigma_v$ be a cohomological irreducible cuspidal automorphic representation of $\GL_3(\A_\F)$ with central character $\omega_\itSigma$. 
We have $|\omega_\itSigma| = |\mbox{ }|_{\A_\F}^{{\sf w}}$ for some even integer $\sf w$.
Let $\mu = \prod_{v \in S_\infty} \in \prod_{v \in S_\infty}X^+(T_3)$ such that $\itSigma$ contributes to the cuspidal cohomology of $\GL_3(\A_\F)$ with coefficients in $M_{\mu,\C}^\vee$.
For each $v \in S_\infty$, we have
\[
\mu_v = (\tfrac{\ell_v-3+{\sf w}}{2},\tfrac{{\sf w}}{2},\tfrac{-\ell_v+3+{\sf w}}{2})
\]
for some odd integer $\ell_v \in \Z_{\geq 3}$ and 
\[
\itSigma_v = {\rm Ind}^{\GL_3(\R)}_{P_{2,1}(\R)}(D_{\ell_v} \boxtimes \epsilon_v) \otimes |\mbox{ }|^{{\sf w}/2}
\]
for some quadratic character $\epsilon_v$ of $\R^\times$. 
Let $L(s,\itSigma \times \itSigma^\vee)$ be the Rankin--Selberg $L$-function of $\itSigma \times \itSigma^\vee$, which has a simple pole at $s=1$. 
We denote by $L^{(\infty)}(s,\itSigma \times \itSigma^\vee)$ the $L$-function obtained by excluding the archimedean $L$-factors.
The aim of this section is to prove Theorem \ref{T:adjoint} on the algebraicity of ${\rm Res}_{s=1}L^{(\infty)}(s,\itSigma \times \itSigma^\vee)$ in terms of product of the Whittaker periods of $\itPi$.
The theorem is a refinement of the result \cite[Theorem 3.3.11]{BR2017} of Balasubramanyam--Raghuram which shall be recalled in Theorem \ref{T:BR} below.
Let $v \in S_\infty$. Let
\[
\<\cdot,\cdot\>_v : \mathcal{W}(\itSigma_v,\psi_{3,\R}) \times \mathcal{W}(\itSigma_v^\vee,\psi_{3,\R}) \longrightarrow \C
\]
be the non-zero $\GL_3(\R)$-equivariant bilinear pairing defined by the local integral
\[
\<W_1,W_2\>_v = \int_{N_2(\R)\backslash \GL_2(\R)} W_1(\iota(g))W_2({\rm diag}(1,-1,1)\iota(g))\,dg,
\]
where $\iota(g) = \bp g & 0 \\ 0 & 1 \ep$. Note that the integral converges absolutely (cf.\,\cite[(3.17)]{JS1981}).
Let 
\[
s: \bigwedge^3 (\frak{g}_{3,\C} / \frak{k}_{3,\C})^* \times \bigwedge^2 (\frak{g}_{3,\C} / \frak{k}_{3,\C})^* \longrightarrow \C
\]
be the ${\rm SO}(3)$-equivariant bilinear pairing defined by
\begin{align*}
X^*\wedge Y^* \wedge Z^*\wedge U^*\wedge V^* &= s(X^*\wedge Y^* \wedge Z^*, U^*\wedge V^*)\\
&\times \bp 1&0&0\\0&0&0\\0&0&0 \ep^*\wedge \bp 1&0&0\\0&1&0\\0&0&0 \ep^* \wedge \bp 0&1&0\\0&0&0\\0&0&0 \ep^*\wedge \bp 0&0&1\\0&0&0\\0&0&0 \ep^*\wedge \bp 0&0&0\\0&0&1\\0&0&0 \ep^*.
\end{align*}
We also fix a non-zero $\GL_3(\C)$-equivariant bilinear pairing 
\[
\<\cdot,\cdot\>_{\mu_v,\C} : M_{\mu_v,\C}^\vee \times M_{\mu_v,\C} \longrightarrow \C
\]
defined over $\Q$.
Let
\[
B(\cdot,\cdot)_v: \left(\mathcal{W}(\itSigma_v,\psi_{3,\R}) \otimes \bigwedge^3 (\frak{g}_{3,\C} / \frak{k}_{3,\C})^* \otimes M_{\mu_v,\C}^\vee\right) \times \left(\mathcal{W}(\itSigma_v^\vee,\psi_{3,\R}) \otimes \bigwedge^2 (\frak{g}_{3,\C} / \frak{k}_{3,\C})^* \otimes M_{\mu_v,\C}\right)\longrightarrow \C
\]
be the bilinear pairing defined by
\[
B(W_1 \otimes X^*\wedge Y^* \wedge Z^*\wedge \otimes P ,\, W_2 \otimes U^*\wedge V^* \otimes Q)_v = \<W_1,W_2\>_v\cdot s(X^*\wedge Y^* \wedge Z^*\wedge U^*\wedge V^*)\cdot \<P,Q\>_{\mu_v,\C}.
\]
For $\sigma \in {\rm Aut}(\C)$, we then have the bilinear pairing
\[
B(\cdot,\cdot)_{\sigma,\infty} = \bigotimes_{v \in S_\infty}B(\cdot,\cdot)_{\sigma^{-1}\circ v}
\]
on 
\[
\left(\mathcal{W}({}^\sigma\!\itSigma_\infty,\psi_{3,\F_\infty}) \otimes \bigwedge^3 (\frak{g}_{3,\C} / \frak{k}_{3,\C})^* \otimes M_{{}^\sigma\!\mu,\C}^\vee\right) \times \left(\mathcal{W}({}^\sigma\!\itSigma_\infty^\vee,\psi_{3,\F_\infty}) \otimes \bigwedge^2 (\frak{g}_{3,\C} / \frak{k}_{3,\C})^* \otimes M_{{}^\sigma\!\mu,\C}\right)
\]
We simply write $B(\cdot,\cdot)_{\sigma,\infty} = B(\cdot,\cdot)_{\infty}$ when $\sigma$ is the identity map.

\begin{thm}[Balasubramanyam--Raghuram]\label{T:BR}
We have
\begin{align*}
&\sigma \left( \frac{{\rm Res}_{s=1}L^{(\infty)}(s,\itSigma \times \itSigma^\vee)}{{\rm Reg}_\F\cdot\pi^d\cdot p^t(\itSigma)\cdot p^b(\itSigma^\vee)}\cdot B([\itSigma_\infty]_t,[\itSigma_\infty^\vee]_b)_{\infty}\right)\\
& = \frac{{\rm Res}_{s=1}L^{(\infty)}(s,{}^\sigma\!\itSigma \times {}^\sigma\!\itSigma^\vee)}{{\rm Reg}_\F\cdot\pi^d\cdot p^t({}^\sigma\!\itSigma)\cdot p^b({}^\sigma\!\itSigma^\vee)}\cdot B([{}^\sigma\!\itSigma_\infty]_t,[{}^\sigma\!\itSigma_\infty^\vee]_b)_{\sigma,\infty}
\end{align*}
for all $\sigma \in {\rm Aut}(\C)$.
\end{thm}

\begin{rmk}
The factor ${\rm Reg}_\F$ appears in the volume of $\F^\times \backslash \A_\F^1$ with respect to the standard measure $\prod_v d^\times x_v$ on $\A_\F^\times$.
As for the factor $\pi^d$, note that we have
\begin{align*}
\int_{\GL_3(\R)}f(g)\,dg = 2^{-1}\Gamma_\R(3)\cdot\int_{\R^2}dx\int_{\R^2}dy\int_{\R^\times}d^\times a \int_{\GL_2(\R)}\frac{dg}{|\det(g)|}\,f\left(\bp {\bf 1}_2 & x \\ 0 & 1\ep \iota(g) \bp {\bf 1}_2 & 0 \\ {}^ty & 1\ep a \right)
\end{align*}
for $f \in L^1(\GL_3(\R))$ (cf.\,\cite[\S\,14.10,\,Corollary]{Voskresenskiibook}).
This formula is used to relate the Rankin--Selberg local zeta integrals for $\itSigma_v \times \itSigma_v^\vee$ at $s=1$ to the bilinear pairing $\<\cdot,\cdot\>_v$ (cf.\,\cite[(3.9)]{Zhang2014} and \cite[(2.2.11)]{BR2017}). 
\end{rmk}

In Theorem \ref{T:adjoint} below, we refine Theorem \ref{T:BR}  by determine the rationality of the archimedean factor $B([\itSigma_\infty]_t,[\itSigma_\infty^\vee]_b)_\infty$ up to an explicit power of $\pi$.
Firstly we give a simplified formula for the archimedean factor in the following lemma.
\begin{lemma}\label{L:archimedean factor 2}
Let $v \in S_\infty$.
We have
\begin{align*}
&B([\itSigma_v]_t,[\itSigma_v^\vee]_b)_v \\
&= (-1)^{{\sf w}/2}4\cdot \frac{(2\ell_v+1)!}{(\ell_v !)^2}\cdot\<W_{(\ell_v, {\sf w}, \epsilon_v; \,0)},W_{(\ell_v, -{\sf w}, \epsilon_v; \,0)}\>_v\cdot \left\<P_{\mu_v^\vee}^+,\,\rho_{\mu_v}\left(\bp 0&0&-1\\0&1&0\\-1&0&0\ep\right) P_{\mu_v}^+ \right\>_{\mu_v,\C}.
\end{align*}
\end{lemma}

\begin{proof}
We drop the subscript $v$ for brevity.
For any ${\rm SO}(3)$-equivariant bilinear pairing  $\<\cdot,\cdot\>:V_\ell \times V_\ell \rightarrow \C$, by (\ref{E:Lie algebra action}), we have
\begin{align}\label{E:Lie algebra bilinear}
\begin{split}
\<E_-^{\ell+i}\cdot {\bf v}_{(\ell;\,\ell)}, E_+^{\ell+i}\cdot {\bf v}_{(\ell;\,-\ell)}\> & = \frac{(2\ell)!(\ell+i)!}{(\ell-i)!}\cdot \<{\bf v}_{(\ell;\,\ell)},{\bf v}_{(\ell;\,-\ell)}\>,\\
\frac{(-1)^{i+1}(2\ell)!}{(\ell+i)!(\ell-i)!}\cdot \<{\bf v}_{(\ell;\,i)},{\bf v}_{(\ell;\,-i)}\>&= \<{\bf v}_{(\ell;\,\ell)},{\bf v}_{(\ell;\,-\ell)}\>
\end{split}
\end{align}
for $-\ell \leq i \leq \ell$.
Let
\[
\<\cdot,\cdot\>_{{\rm ad},\mu} : \left(\bigwedge^3 (\frak{g}_{3,\C} / \frak{k}_{3,\C})^* \otimes M_{\mu,\C}^\vee\right) \times \left(\bigwedge^2 (\frak{g}_{3,\C} / \frak{k}_{3,\C})^* \otimes M_{\mu,\C}\right) \longrightarrow \C
\]
be the ${\rm SO}(3)$-equivariant bilinear pairing defined by
\[
\<X^*\wedge Y^* \wedge Z^* \otimes P,\,U^*\wedge V^* \otimes Q\>_{{\rm ad},\mu} = s(X^*\wedge Y^* \wedge Z^*, U^*\wedge V^*) \cdot \<P,Q\>_{\mu,\C}.
\]
As explained in the proof of Lemma \ref{L:cohomology generator}, we have ${\rm SO}(3)$-equivariant homomorphisms 
\begin{align*}
V_\ell &\longrightarrow \bigwedge^2 (\frak{g}_{3,\C} / \frak{k}_{3,\C})^* \otimes M_{\mu,\C}, \quad {\bf v}_{(\ell;\,-\ell)} \longmapsto X_{1}^* \wedge X_{2}^* \otimes \rho_{\mu}\left( \bp -1 & 0 & -1  \\ \sqrt{-1} & 0 & -\sqrt{-1} \\ 0 & 1 & 0\ep\right)P_{\mu}^+,\\
V_\ell &\longrightarrow \bigwedge^3 (\frak{g}_{3,\C} / \frak{k}_{3,\C})^* \otimes M_{\mu,\C}^\vee, \quad {\bf v}_{(\ell;\,\ell)} \longmapsto X_0^*\wedge X_{-1}^* \wedge X_{-2}^* \otimes \rho_{\mu^\vee}\left( \bp 1 & 0 & 1  \\ \sqrt{-1} & 0 & -\sqrt{-1} \\ 0 & 1 & 0\ep\right)P_{\mu^\vee}^+.
\end{align*}
Note that $s(X_0^*\wedge X_{-1}^* \wedge X_{-2}^*,X_1^*\wedge X_2^*) = -4\sqrt{-1}$ and
\begin{align*}
&\left\<\rho_{\mu^\vee}\left( \bp 1 & 0 & 1  \\ \sqrt{-1} & 0 & -\sqrt{-1} \\ 0 & 1 & 0\ep\right)P_{\mu^\vee}^+,\, \rho_{\mu}\left( \bp -1 & 0 & -1  \\ \sqrt{-1} & 0 & -\sqrt{-1} \\ 0 & 1 & 0\ep\right)P_{\mu}^+\right\>_{\mu,\C}\\
& = \left\<P_{\mu^\vee}^+,\, \rho_{\mu}\left( \bp 0 & 0 & -1  \\ 0 & 1 & 0 \\ -1 & 0 & 0\ep\right)P_{\mu}^+\right\>_{\mu,\C} \in \Q^\times.
\end{align*}
It then follows from (\ref{E:relation}) and (\ref{E:Lie algebra bilinear}) that
\begin{align*}
\<[\itSigma_\infty]_t,[\itSigma_\infty^\vee]_b\> &= 4\cdot(2\ell)!\cdot\left\<P_{\mu^\vee}^+,\, \rho_{\mu}\left( \bp 0 & 0 & -1  \\ 0 & 1 & 0 \\ -1 & 0 & 0\ep\right)P_{\mu}^+\right\>_{\mu,\C}\\
&\times\sum_{i=-\ell}^{\ell}\frac{(-1)^{i+{\sf w}/2}}{(\ell+i)!(\ell-i)!} \<W_{(\ell, {\sf w}, \epsilon; \,i)},W_{(\ell, -{\sf w}, \epsilon; \,-i)}\>_\infty\\
&=(-1)^{{\sf w}/2}4\cdot \frac{(2\ell+1)!}{(\ell !)^2}\cdot\<W_{(\ell, {\sf w}, \epsilon; \,0)},W_{(\ell, -{\sf w}, \epsilon; \,0)}\>_\infty\cdot \left\<P_{\mu^\vee}^+,\,\rho_\mu\left(\bp 0&0&-1\\0&1&0\\-1&0&0\ep\right) P_{\mu}^+ \right\>_{\mu,\C}.
\end{align*}
This completes the proof.
\end{proof}

In the following lemma, we compute the archimedean pairing of Whittaker functions that appearing in the simplified formula for $B([\itSigma_\infty]_t,[\itSigma_\infty^\vee]_b)_\infty$.
\begin{lemma}\label{L:archimedean adjoint}
Let $v \in S_\infty$.
We have
\[
\<W_{(\ell_v, {\sf w}, \epsilon_v; \,0)},W_{(\ell_v, -{\sf w}, \epsilon_v; \,0)}\>_v = -4\cdot \frac{\Gamma_\C(\ell_v+1)}{\Gamma_\C(1)\Gamma_\R(2\ell_v+3)}\cdot L(1,\itSigma_v \times \itSigma_v^\vee).
\]
\end{lemma}

\begin{proof}
We may assume ${\sf w}=0$ after replacing $\itSigma_\infty$ by $\itSigma_\infty \otimes |\mbox{ }|^{-{\sf w}/2}$. Write $W_{(\ell;\, 0)} = W_{(\ell, {\sf w}, \epsilon;\, 0)}$.
Let $L_1$ and $L_2$ be vertical paths from $c_1-\sqrt{-1}\,\infty$ to $c_1+\sqrt{-1}\,\infty$ and $c_2-\sqrt{-1}\,\infty$ to $c_2+\sqrt{-1}\,\infty$, respectively, for some
\[
0<c_1<1,\quad -(\tfrac{\ell-1}{2})<c_2<\tfrac{\ell+1}{2}.
\]
We have
\begin{align*}
& \<W_{(\ell;\, 0)},W_{(\ell;\, 0)}\>_\infty\\
& = \int_{N_2(\R)\backslash \GL_2(\R)} W_{(\ell;\, 0)}(\iota(g))W_{(\ell;\, 0)}({\rm diag}(1,-1,1)\iota(g))\,dg\\
& = \int_0^\infty \frac{d^\times a_1}{a_1}\int_0^\infty d^\times a_2\,W_{(\ell;\, 0)}\left(\iota\left({\rm diag}(a_1a_2,a_2) \right)\right)W_{(\ell;\, 0)}\left(\iota\left({\rm diag}(a_1a_2,a_2) \right)\right)\\
& = -\int_0^\infty d^\times a_1\int_0^\infty d^\times a_2\int_{L_1}\frac{ds_1}{2\pi\sqrt{-1}}\int_{L_2}\frac{ds_2}{2\pi\sqrt{-1}}\int_{L_1}\frac{ds_1'}{2\pi\sqrt{-1}}\int_{L_2}\frac{ds_2'}{2\pi\sqrt{-1}}\,a_1^{-s_1-s_1'+1}a_2^{-s_2-s_2'+2}\\
& \quad\quad\quad\quad\quad\quad\times \frac{\Gamma_\C(s_1+\tfrac{\ell-1}{2})\Gamma_\R(s_1)\Gamma_\C(s_2+\tfrac{\ell-1}{2})\Gamma_\R(s_2+\ell)\Gamma_\C(s_1'+\tfrac{\ell-1}{2})\Gamma_\R(s_1')\Gamma_\C(s_2'+\tfrac{\ell-1}{2})\Gamma_\R(s_2'+\ell)}{\Gamma_\R(s_1+s_2+\ell)\Gamma_\R(s_1'+s_2'+\ell)}\\
& = -\int_{L_2}\frac{ds_2}{2\pi\sqrt{-1}}\Gamma_\C(s_2+\tfrac{\ell-1}{2})\Gamma_\R(s_2+\ell)\Gamma_\C(-s_2+\tfrac{\ell+3}{2})\Gamma_\R(-s_2+\ell+2)\\
&\quad\times\int_{L_1}\frac{ds_1}{2\pi\sqrt{-1}}\frac{\Gamma_\C(s_1+\tfrac{\ell-1}{2})\Gamma_\R(s_1)\Gamma_\C(-s_1+\tfrac{\ell+1}{2})\Gamma_\R(-s_1+1)}{\Gamma_\R(s_1+s_2+\ell)\Gamma_\R(-s_1-s_2+\ell+3)}.
\end{align*}
Here the last equality follows from the Mellin inversion formula. By the first Barnes lemma, we have
\begin{align*}
&\int_{L_1'}\frac{dt_1}{2\pi\sqrt{-1}}\,\Gamma_\R(t_1+s_1)\Gamma_\R(t_1+s_2+\tfrac{\ell+1}{2})\Gamma_\R(-t_1+\tfrac{\ell-1}{2})\Gamma_\R(-t_1)\\
& = 2\cdot \frac{\Gamma_\R(s_1+\tfrac{\ell-1}{2})\Gamma_\R(s_1)\Gamma_\R(s_2+\ell)\Gamma_\R(s_2+\tfrac{\ell+1}{2})}{\Gamma_\R(s_1+s_2+\ell)},\\
&\int_{L_2'}\frac{dt_2}{2\pi\sqrt{-1}}\,\Gamma_\R(t_2-s_1)\Gamma_\R(t_2-s_2+\tfrac{\ell+3}{2})\Gamma_\R(-t_2+\tfrac{\ell+1}{2})\Gamma_\R(-t_2+1)\\
 &= 2\cdot\frac{\Gamma_\R(-s_1+\tfrac{\ell+1}{2})\Gamma_\R(-s_1+1)\Gamma_\R(-s_2+\ell+2)\Gamma_\R(-s_2+\tfrac{\ell+5}{2})}{\Gamma_\R(-s_1-s_2+\ell+3)}.
\end{align*}
Here $L_1'$ and $L_2'$ are vertical paths from $c_1'-\sqrt{-1}\,\infty$ to $c_1'+\sqrt{-1}\,\infty$ and $c_2'-\sqrt{-1}\,\infty$ to $c_2'+\sqrt{-1}\,\infty$, respectively, for some 
\[
\max\{-c_1,-c_2-(\tfrac{\ell+1}{2})\} < c_1' <0,\quad \max\{c_1,c_2-(\tfrac{\ell+3}{2})\} < c_2' < 1.
\]
Therefore we have
\begin{align*}
&\int_{L_1}\frac{ds_1}{2\pi\sqrt{-1}}\frac{\Gamma_\C(s_1+\tfrac{\ell-1}{2})\Gamma_\R(s_1)\Gamma_\C(-s_1+\tfrac{\ell+1}{2})\Gamma_\R(-s_1+1)}{\Gamma_\R(s_1+s_2+\ell)\Gamma_\R(-s_1-s_2+\ell+3)}\\
& = 2^{-2}\cdot \frac{1}{\Gamma_\R(s_2+\ell)\Gamma_\R(s_2+\tfrac{\ell+1}{2})\Gamma_\R(-s_2+\ell+2)\Gamma_\R(-s_2+\tfrac{\ell+5}{2})}\\
& \times \int_{L_1'}\frac{dt_1}{2\pi\sqrt{-1}}\,\Gamma_\R(t_1+s_2+\tfrac{\ell+1}{2})\Gamma_\R(-t_1+\tfrac{\ell-1}{2})\Gamma_\R(-t_1)\\
&\times \int_{L_2'}\frac{dt_2}{2\pi\sqrt{-1}}\,\Gamma_\R(t_2-s_2+\tfrac{\ell+3}{2})\Gamma_\R(-t_2+\tfrac{\ell+1}{2})\Gamma_\R(-t_2+1)\\
&\times \int_{L_1}\frac{ds_1}{2\pi\sqrt{-1}}\,\Gamma_\R(s_1+t_1)\Gamma_\R(s_1+\tfrac{\ell+1}{2})\Gamma_\R(-s_1+t_2)\Gamma_\R(-s_1+\tfrac{\ell+3}{2})\\
& = 2^{-1}\cdot \frac{\Gamma_\R(\ell+2)}{\Gamma_\R(s_2+\ell)\Gamma_\R(s_2+\tfrac{\ell+1}{2})\Gamma_\R(-s_2+\ell+2)\Gamma_\R(-s_2+\tfrac{\ell+5}{2})}\int_{L_1'}\frac{dt_1}{2\pi\sqrt{-1}}\int_{L_2'}\frac{dt_2}{2\pi\sqrt{-1}}\,\frac{\Gamma_\R(t_1+t_2)}{\Gamma_\R(t_1+t_2+\ell+2)}\\
& \quad\quad\quad\quad\quad\quad\quad\quad\quad\quad\quad\quad\quad\quad\quad\quad\quad\quad\quad\quad\times \Gamma_\R(t_1+s_2+\tfrac{\ell+1}{2})\Gamma_\R(t_1+\tfrac{\ell+3}{2})\Gamma_\R(-t_1+\tfrac{\ell-1}{2})\Gamma_\R(-t_1)\\
&\quad\quad\quad\quad\quad\quad\quad\quad\quad\quad\quad\quad\quad\quad\quad\quad\quad\quad\quad\quad\times \Gamma_\R(t_2-s_2+\tfrac{\ell+3}{2})\Gamma_\R(t_2+\tfrac{\ell+1}{2})\Gamma_\R(-t_2+\tfrac{\ell+1}{2})\Gamma_\R(-t_2+1).
\end{align*}
Here the last equality follows from the first Barnes lemma.
Hence $\<W_{(\ell;\, 0)},W_{(\ell;\, 0)}\>_\infty$ is equal to
\begin{align*}
& -2^{-1}\cdot\Gamma_\R(\ell+2)\int_{L_1'}\frac{dt_1}{2\pi\sqrt{-1}}\int_{L_2'}\frac{dt_2}{2\pi\sqrt{-1}}\,\frac{\Gamma_\R(t_1+t_2)}{\Gamma_\R(t_1+t_2+\ell+2)}\\
&\quad\quad\quad\quad\quad\quad\quad\quad\times \Gamma_\R(t_1+\tfrac{\ell+3}{2})\Gamma_\R(-t_1+\tfrac{\ell-1}{2})\Gamma_\R(-t_1)\Gamma_\R(t_2+\tfrac{\ell+1}{2})\Gamma_\R(-t_2+\tfrac{\ell+1}{2})\Gamma_\R(-t_2+1)\\
&\quad\quad\quad\quad\quad\quad\quad\quad\times \int_{L_2}\frac{ds_2}{2\pi\sqrt{-1}}\,\Gamma_\R(s_2+t_1+\tfrac{\ell+1}{2})\Gamma_\R(s_2+\tfrac{\ell-1}{2})\Gamma_\R(-s_2+t_2+\tfrac{\ell+3}{2})\Gamma_\R(-s_2+\tfrac{\ell+3}{2})\\
& =-\Gamma_\R(\ell+1)\Gamma_\R(\ell+2)\int_{L_1'}\frac{dt_1}{2\pi\sqrt{-1}}\,\Gamma_\R(t_1+\tfrac{\ell+3}{2})\Gamma_\R(t_1+\ell+2)\Gamma_\R(-t_1+\tfrac{\ell-1}{2})\Gamma_\R(-t_1)\\
&\quad\quad\quad\quad\quad\quad\quad\quad\quad\times \int_{L_2'}\frac{dt_2}{2\pi\sqrt{-1}}\,\frac{\Gamma_\R(t_2+t_1)\Gamma_\R(t_2+\ell+1)\Gamma_\R(t_2+\tfrac{\ell+1}{2})\Gamma_\R(-t_2+\tfrac{\ell+1}{2})\Gamma_\R(-t_2+1)}{\Gamma_\R(t_2+t_1+2\ell+3)}\\
& = -2\cdot\frac{\Gamma_\R(\tfrac{\ell+3}{2})\Gamma_\R(\tfrac{3\ell+3}{2})\Gamma_\R(\ell+1)^2\Gamma_\R(\ell+2)^2}{\Gamma_\R(2\ell+3)}\\
& \times \int_{L_1'}\frac{dt_1}{2\pi\sqrt{-1}}\,\frac{\Gamma_\R(t_1+\tfrac{\ell+1}{2})\Gamma_\R(t_1+\tfrac{\ell+3}{2})\Gamma_\R(t_1+1)\Gamma_\R(-t_1+\tfrac{\ell-1}{2})\Gamma_\R(-t_1)}{\Gamma_\R(t_1+\tfrac{3\ell+5}{2})}\\
& = -4\cdot \frac{\Gamma_\R(1)\Gamma_\C(\ell)\Gamma_\C(\ell+1)\Gamma_\C(\tfrac{\ell+1}{2})^2}{\Gamma_\R(2\ell+3)}.
\end{align*}
Here the first equality  and the last two equalities follow from Barnes' first and second lemmas, respectively.
Finally, note that
\[
L(s,\itSigma_v \times \itSigma_v^\vee) = \Gamma_\R(s)\Gamma_\C(s)\Gamma_\C(s+\ell-1)\Gamma_\C(s+\tfrac{\ell-1}{2})^2.
\]
This completes the proof.
\end{proof}

\begin{thm}\label{T:adjoint}
We have
\[
\sigma \left(\frac{{\rm Res}_{s=1}L^{(\infty)}(s,\itSigma \times \itSigma^\vee)}{{\rm Reg}_\F\cdot \pi^{2d+2\sum_{v \in S_\infty}\ell_v}\cdot p^t(\itSigma)\cdot p^b(\itSigma^\vee)}\right) = \frac{{\rm Res}_{s=1}L^{(\infty)}(s,{}^\sigma\!\itSigma \times {}^\sigma\!\itSigma^\vee)}{{\rm Reg}_\F\cdot \pi^{2d+2\sum_{v \in S_\infty}\ell_v}\cdot p^t(\itSigma)\cdot p^b({}^\sigma\!\itSigma^\vee)} 
\]
for all $\sigma \in {\rm Aut}(\C)$.
\end{thm}

\begin{proof}
By the definition of the classes $[\itSigma_\infty]_b$, $[{}^\sigma\!\itSigma_\infty]_b$, $[\itSigma_\infty^\vee]_t$, and $[{}^\sigma\!\itSigma_\infty^\vee]_t$, it is clear that
\[
B([{}^\sigma\!\itSigma_\infty]_t,[{}^\sigma\!\itSigma_\infty^\vee]_b)_{\sigma,\infty} = B([\itSigma_\infty]_t,[\itSigma_\infty^\vee]_b)_{\infty}.
\]
By Lemmas \ref{L:archimedean factor 2} and \ref{L:archimedean adjoint}, we have
\[
B([\itSigma_v]_t,[\itSigma_v^\vee]_b)_v \in \pi^{-2\ell_v-1}\cdot \Q^\times
\]
for all $v \in S_\infty$.
The assertion then follows from Theorem \ref{T:BR}.
This completes the proof.
\end{proof}


\end{document}